\documentclass[12pt]{article}

\usepackage{amssymb}

\newtheorem {theorem} {Theorem}%[section]
\newtheorem {definition} [theorem]{Definition}
\newtheorem {proposition} [theorem]{Proposition}

\newtheorem {lemma}  [theorem]{Lemma}
\newtheorem {remark} [theorem]{\sc Remark}

\newcommand{\bbox}{\ \hfill\rule[-1mm]{2mm}{3.2mm}}

\title{The role of algebraic solutions in planar polynomial differential
systems. \thanks{The second and third authors are partially
supported by a MCYT grant number BFM 2002-04236-C02-01. The second
author is also partially supported by DURSI of Government of
Catalonia ``Distinci\'o de la Generalitat de Catalunya per a la
promoci\'o de la recerca universit\`aria''.}}

\author{{\sc H\'ector Giacomini$^{\ (1)}$, Jaume Gin\'e$^{\ (2)}$ and Maite Grau$^{\ (2)}$}}
\date{}

\begin{document}
\maketitle

\begin{abstract}
We study a planar polynomial differential system, given by
$\dot{x}=P(x,y)$, $\dot{y}=Q(x,y)$. We consider a function
$I(x,y)=\exp \{ h_2(x) A_1(x,y) \diagup A_0(x,y) \}$ $ h_1(x)
\prod_{i=1}^{\ell} (y-g_i(x))^{\alpha_i}$, where $g_i(x)$ are
algebraic functions of $x$, $A_1(x,y)=\prod_{k=1}^r (y-a_k(x))$,
$A_0(x,y)=\prod_{j=1}^s (y-\tilde{g}_j(x))$ with $a_k(x)$ and
$\tilde{g}_j(x)$ algebraic functions, $A_0(x,y)$ and $A_1(x,y)$ do
not share any common factor, $h_2(x)$ is a rational function,
$h(x)$ and $h_1(x)$ are functions of $x$ with a rational
logarithmic derivative and $\alpha_i \in \mathbb{C}$. We show that
if $I(x,y)$ is a first integral or an integrating factor, then
$I(x,y)$ is a Darboux function. A Darboux function is a function
of the form $f_1^{\lambda_1} \ldots f_p^{\lambda_p} \exp\{
h/f_0\}$, where $f_i$ and $h$ are polynomials in $\mathbb{C}[x,y]$
and the $\lambda_i$'s are complex numbers. In order to prove this
result, we show that if $g(x)$ is an algebraic particular
solution, that is, if there exists an irreducible polynomial
$f(x,y)$ such that $f(x,g(x)) \equiv 0$, then $f(x,y)=0$ is an
invariant algebraic curve of the system. In relation with this
fact, we give some characteristics related to particular solutions
and functions of the form $I(x,y)$ such as the structure of their
cofactor.
\par
Moreover, we consider $A_0(x,y)$, $A_1(x,y)$ and $h_2(x)$ as
before and a function of the form $\Phi(x,y):= \exp \{h_2(x)\,
A_1(x,y) / A_0 (x,y) \}$. We show that if the derivative of
$\Phi(x,y)$ with respect to the flow is well defined over $\{
(x,y) \, :\, A_0(x,y)=0 \}$ then $\Phi(x,y)$ gives rise to an
exponential factor. This exponential factor has the form
$\exp\{R(x,y)\}$ where $R=h_2 A_1/A_0 + B_1/B_0$ and with
$B_1/B_0$ a function of the same form as $h_2 A_1/A_0$. Hence,
$\exp\{ R(x,y)\}$ factorizes as the product $\Phi(x,y) \,
\Psi(x,y)$, for $\Psi(x,y):=\exp\{ B_1(x,y) / B_0(x,y)\}$.
\end{abstract}

{\small{\noindent 2000 {\it AMS Subject Classification:} 14H05, 34A05, 34A34  \\
\noindent {\it Key words and phrases:} Planar polynomial
differential system, algebraic function, invariant algebraic
curve, integrability.}}

\section{Introduction \label{sect1}}

In this work we consider planar polynomial differential systems
as:
\begin{equation}
\dot{x}= P(x,y), \quad  \dot{y}=Q(x,y), \label{eq1}
\end{equation}
where $P(x,y)$ and $Q(x,y)$ belong to the ring of real polynomials
in two variables, $\mathbb{R}[x,y]$. We will always assume that
$P(x,y)$ and $Q(x,y)$ are coprime polynomials. We denote by ${\rm
d}= \max \{ \deg P, \deg Q \}$ and we say that ${\rm d}$ is the
degree of system (\ref{eq1}). Equivalently to system (\ref{eq1}),
we may consider the ordinary differential equation:
\begin{equation}
\frac{dy}{dx} \, = \, \frac{Q(x,y)}{P(x,y)}. \label{eq2}
\end{equation}
In order to shorten some formulae, we introduce the operator
$\mathcal{X}$ associated to (\ref{eq1}):
\[ \mathcal{X} := P(x,y) \frac{\partial }{\partial x} + Q(x,y) \frac{\partial }{\partial
y}. \]
\par
This paper is related with the study of the properties of a
certain type of particular solutions of system (\ref{eq1}). Given
a function $g(x)$, we say that $y=g(x)$ is a {\em particular
solution} of equation (\ref{eq2}) if
\begin{equation}
g'(x) \, = \, \frac{Q(x,g(x))}{P(x,g(x))}, \label{eq3}
\end{equation}
where $g'(x)=dg(x)/dx$. In particular, we are concerned with {\em
algebraic} functions. We say that a function $g(x)$ is algebraic
if there exists a polynomial $f(x,y)$ such that $f(x,g(x)) \equiv
0$. As we prove in the next section, by using some algebraic
results stated for instance in \cite{Casas, Hille, Walker}, this
polynomial $f(x,y)$ can always be chosen irreducible and it is
unique modulus multiplication by constants.

The solutions of system (\ref{eq1}) may also be given in an
implicit way. An {\em invariant algebraic curve} of system
(\ref{eq1}) is an algebraic curve $f(x,y)=0$ satisfying
\begin{equation} \left. P(x,y) \frac{\partial f}{\partial x} (x,y) +
Q(x,y) \frac{\partial f}{\partial y}(x,y) \right|_{f(x,y)=0} =0.
\label{eq3a}
\end{equation}
By using Hilbert's Nullstellensatz, cf. \cite{Lang}, it can be
shown that $f(x,y)=0$ is an invariant algebraic curve of system
(\ref{eq1}) if, and only if, there exists a polynomial $k(x,y) \in
\mathbb{R}[x,y]$ satisfying:
\[ P(x,y) \frac{\partial f}{\partial x} (x,y) + Q(x,y) \frac{\partial
f}{\partial y}(x,y) = k(x,y) f(x,y). \] This polynomial $k(x,y)$
is called the {\em cofactor} of the curve given by $f(x,y)=0$ and
it can be shown that its degree is lower than or equal to ${\rm
d}-1$. Invariant algebraic curves, also denoted as Darboux
polynomials by some authors, have been widely studied for their
relation with integrability and some qualitative properties of
polynomial differential systems, see for instance \cite{Moulin1,
Moulin2}. As a generalization of the notion of invariant algebraic
curve, we can define an {\em exponential factor}. This concept has
been firstly introduced by Christopher \cite{Christopher1} and it
is related with the notion of multiplicity of an invariant
algebraic curve of system (\ref{eq1}). An exponential factor is a
function of the form $\exp \{ h/f_0 \}$ where $h(x,y)$ and
$f_0(x,y) \in \mathbb{C}[x,y]$, and $\mathcal{X} \left( \exp \{h /
f_0 \} \right) \, = \, \tilde{k}(x,y) \, \exp \{h / f_0 \}$ with
$\tilde{k}(x,y)$ a polynomial of degree at most ${\rm d}-1$. The
following lemma is given in \cite{Christopher1} and characterizes
exponential factors.
\begin{lemma}
The function $\exp \{h / f_0 \}$ is an exponential factor of
system {\rm (\ref{eq1})} with cofactor $\tilde{k}(x,y)$ if, and
only if, $f_0(x,y)=0$ is an invariant algebraic curve of system
{\rm (\ref{eq1})} with cofactor $k_0(x,y)$ and $\mathcal{X}(h) =
k_0 \, h \, + \, \tilde{k} \, f_0$. \label{lem1}
\end{lemma}
Invariant algebraic curves characterize the existence of first
integrals for system (\ref{eq1}) belonging to a certain functional
class. In order to properly state the known results about
integrability using invariant algebraic curves, we need to
consider complex algebraic curves $f(x,y)=0$, where $f(x,y) \in
\mathbb{C}[x,y]$. Since system (\ref{eq1}) is defined by real
polynomials, if $f(x,y)=0$ is an invariant algebraic curve with
cofactor $k(x,y)$, then its conjugate $\bar{f}(x,y)=0$ is also an
invariant algebraic curve with cofactor $\bar{k}(x,y)$. Hence, its
product $f(x,y) \bar{f}(x,y) \in \mathbb{R}[x,y]$ gives rise to a
real invariant algebraic curve with a real cofactor $k(x,y)+
\bar{k}(x,y)$. In the Darboux theory of integrability, quoted in
the forthcoming paragraph, we need to consider invariant algebraic
curves defined by polynomials in $\mathbb{C}[x,y]$ since they play
an essential role in the theory of integrability. We also notice
that in $\mathbb{R}^2$ the curve given by $f(x,y)=0$, even if
$f(x,y)\in \mathbb{R}[x,y]$, may only contain a finite number of
isolated singular points or be the null set.
\par
A function of the form $f_1^{\lambda_1} \ldots f_p^{\lambda_p}
\exp\{ h/f_0\}$, where $f_i$ and $h$ are polynomials in
$\mathbb{C}[x,y]$ and the $\lambda_i \in \mathbb{C}$, is called a
{\it Darboux function}, see for instance \cite{Moulin1, Moulin2}.
System (\ref{eq1}) is called {\it Darboux integrable} if it has a
first integral or an integrating factor which is a Darboux
function (for a definition of first integral and of integrating
factor, see \cite{ChGGLl}). The following lemma gives the relation
between Darboux functions and invariant algebraic curves of a
system (\ref{eq1}).
\begin{lemma}
We consider a Darboux function $\mathcal{D}(x,y):=f_1^{\lambda_1}
\ldots f_p^{\lambda_p} \exp\{h/f_0\}$ such that
\[ P(x,y) \frac{\partial \mathcal{D}}{\partial x}(x,y) + Q(x,y) \frac{\partial
\mathcal{D}}{\partial y}(x,y) = k(x,y) \mathcal{D}(x,y), \] where
$k(x,y)$ is a polynomial of degree at most ${\rm d}-1$. Then, each
$f_i(x,y)=0$, $i=1,2, \ldots, p$ is an invariant algebraic curve
of system {\rm (\ref{eq1})} and $\exp\{h/f_0\}$ is an exponential
factor of system {\rm (\ref{eq1})}. \label{lem2} \
\end{lemma}
The proof of this lemma is analogous to the proofs of Lemma 3 and
Proposition 4 in \cite{ChLl}. \newline

We recall that $V(x,y)$ is an inverse integrating factor of system
(\ref{eq1}) if it is a function of class $\mathcal{C}^1$ in some
open set $\mathcal{U}$ of $\mathbb{R}^2$ and satisfies the
following partial differential equation:
\[ P(x,y) \frac{\partial V}{\partial x}(x,y) + Q(x,y) \frac{\partial
V}{\partial y}(x,y) = \left( \frac{\partial P}{\partial x}(x,y) +
\frac{\partial Q}{\partial y}(x,y) \right) \,V(x,y) . \] We note
that the function $1/V(x,y)$ is an integrating factor for system
(\ref{eq1}) in $\mathcal{U}$. The following result, which is a
summary of some well known results, relates the existence of an
inverse integrating factor in a certain functional class and the
existence of a first integral in another (possibly larger)
functional class. The definitions of elementary function and
Liouvillian function can be found in \cite{Singer}. Proposition
\ref{prop3} shows that when considering the integrability problem
we are also addressed to study whether an inverse integrating
factor belongs to a certain given class of functions. As many
authors have noted, see for instance \cite{ChGGLl}, inverse
integrating factors play a fundamental role in the integrability
problem, not only because they characterize the functional class
of a first integral but also because they usually belong to an
easier functional class. For instance, quadratic systems of the
form (\ref{eq1}) with a center at the origin always have an
inverse integrating factor which is a polynomial of degree $3$ or
$5$. Hence, we have that the characterization of centers for
quadratic systems can be done by means of a polynomial instead of
a first integral, which will be of Darboux type in a general case.
\begin{proposition} The following three statements hold.
\begin{itemize}

\item[{\rm (a)}] System {\rm (\ref{eq1})} has a {\em Darboux}
first integral if, and only if, it has a {\em rational} inverse
integrating factor.

\item[{\rm (b)}] If system {\rm (\ref{eq1})} has an {\em
elementary} first integral, then it has an inverse integrating
factor of the form $V(x,y)=\left( A(x,y) /B(x,y) \right)^{1/N}$,
where $N \in \mathbb{Z}$ and $A,B \in \mathbb{C}[x,y]$.

\item[{\rm (c)}] System {\rm (\ref{eq1})} has a {\em Liouvillian}
first integral if, and only if, it has a {\em Darboux} inverse
integrating factor.
\end{itemize} \label{prop3}
\end{proposition}
The first statement of this proposition is proved in \cite{ChGGLl}
and its reciprocal is proved in \cite{Christopher2, Schinzel}.
Statement (b) is proved in \cite{PrelleSinger} and the last
statement is proved in \cite{Singer}.
\newline

It is clear, as shown in Lemma $2$ of \cite{GaGG}, that given an
invariant algebraic curve $f(x,y)=0$ of system (\ref{eq1}), all
the algebraic functions defined by it in an implicit way, that is,
all the functions $g(x)$ satisfying $f(x,g(x)) \equiv 0$, are
particular solutions of equation (\ref{eq2}). In this paper, among
other results, the converse result is established, that is, given
an algebraic particular solution $y=g(x)$ of equation (\ref{eq2}),
we show that the irreducible polynomial $f(x,y)$ such that
$f(x,g(x)) \equiv 0$ gives rise to an invariant algebraic curve
$f(x,y)=0$ of system (\ref{eq1}). This fact is stated and proved
in Theorem \ref{th6} of Section \ref{sect3}. We have noticed that
a Darboux function may also contain exponential factors and this
fact is necessary so as to characterize the Liouvillian
integrability. Hence, exponential factors appear in a natural way
when considering invariant algebraic curves. In this paper we
consider algebraic particular solutions $y-g(x)$ which come
naturally from invariant algebraic curves and this relationship
allows us to give an analogous to exponential factors but for
algebraic particular solutions. This analogy is motivated and made
clear in Subsection \ref{sect32}. In relation with this fact, we
give some characteristics related to particular solutions such as
the structure of their cofactor, which is given in Subsection
\ref{sect33}.
\newline

In \cite{GaGG} and \cite{GG} an algorithmic method to determine,
for system (\ref{eq1}), the possible existence of first integrals
or integrating factors of the form $I(x,y)= h(x)
\prod_{i=1}^{\ell} (y-g_i(x))^{\alpha_i}$ or $I(x,y)= \exp \{
h_2(x) \prod_{k=1}^r (y-a_k(x)) \diagup
\prod_{j=1}^s(y-\tilde{g}_j(x)) \} $ $h_1(x)$ $\prod_{i=1}^{\ell}
(y-g_i(x))^{\alpha_i}$, where $g_i(x)$ and $\tilde{g}_j(x)$ are
unknown particular solutions of equation (\ref{eq2}), $\alpha_i
\in \mathbb{C}$ are unknown constants, and $a_k(x)$, $h(x)$,
$h_1(x)$ and $h_2(x)$ are unknown functions, is given. In both
cases, if all the particular solutions $g_i(x)$ and
$\tilde{g}_j(x)$ are determined, which is expressed by the
non-existence of a nonlinear superposition principle as described
in \cite{GaGG, GG}, they are algebraic functions (see Proposition
7 in \cite{GaGG} and Theorem 2 in \cite{GG}). The algorithm, in
this case, gives an alternative method to determine such type of
solutions. In the case where all the $g_i(x)$, $a_k(x)$ and
$\tilde{g}_j(x)$ are algebraic functions, $h(x)$ and $h_1(x)$ have
a rational logarithmic derivative and $h_2(x)$ is a rational
function, we show (cf. Propositions \ref{prop16} and \ref{prop17},
Subsection \ref{sect34}) that $I(x,y)$ is a Darboux function. This
result is one of the main goals obtained in the present work.
Hence, using the algorithm described in \cite{GaGG, GG}, all the
systems with a Liouvillian first integral can be found, as well as
an explicit expression of a non-Liouvillian first integral when
there is a nonlinear superposition principle. The present work is
born as a complement to the results described in \cite{GaGG, GG}.
In these two works the integrability problem is studied in the
particular case of the existence of a first integral of the
described form $I(x,y)$. As we have already stated, in \cite{GaGG,
GG} it is shown that if when applying the described algorithm the
function $I(x,y)$ is completely determined then it only involves
algebraic functions. In that case, in all the examples studied in
\cite{GaGG, GG}, the function $I(x,y)$ was a Darboux function. In
fact, this is the general case: if when applying the algorithm
described in \cite{GaGG, GG} the function $I(x,y)$ is completely
determined then it is a Darboux function. The proof of this
assertion is one of the objectives that we have achieved in the
present work.

\section{Some preliminary results on algebraic functions and polynomials \label{sect2}}

In this section we give a summary of well known algebraic results
which are needed in this paper. The definitions and proofs can be
found, for instance, in the books \cite{Casas, Lang, Walker}.

We always consider polynomials in $\mathbb{R}[x,y]$ which is a
ring with the usual addition and product operations of
polynomials. Equivalently, we may consider the ring
$\mathbb{C}[x,y]$. One of the most important equivalence relations
which can be defined in the ring $\mathbb{R}[x,y]$ is the
divisibility relation. We say that the polynomial $f_1 \in
\mathbb{R}[x,y]$ divides $f_2 \in \mathbb{R}[x,y]$, and we write
$f_1 \, |\ f_2$, if there exists a polynomial $k \in
\mathbb{R}[x,y]$ such that $f_2=k \, f_1$. It can be shown that
$\mathbb{R}[x,y]$ is a unique factorization domain (UFD). We
recall that the unit elements of the ring $\mathbb{R}[x,y]$ are
the constants different from zero, i.e., $\mathbb{R}- \{ 0 \}$. We
say that two elements $f_1,f_2 \in \mathbb{R}[x,y]$ are {\em
associates} if there exists a unit element $e$ such that $f_1= e
\, f_2$. An irreducible polynomial is a non--constant element
$f(x,y) \in \mathbb{R}[x,y]$ which is only divided by its
associates.
\newline

If $f_1, f_2 \in \mathbb{R}[x,y]$ are such that for any $(a,b) \in
\{ (x,y) \in \mathbb{C}^2 \, : \ f_1(x,y)=0 \}$ we have that
$f_2(a,b)=0$, we will write $\left. f_2(x,y)
\right|_{f_1(x,y)=0}=0$. Given a polynomial $f \in
\mathbb{R}[x,y]$, it defines a curve, denoted by $f(x,y)=0$, which
is the set $\{ (x,y) \in \mathbb{C}^2 \, : \ f(x,y)=0 \}$. We say
that $p \in \mathbb{C}^2$ is a point of intersection of two curves
$f_1(x,y)=0$ and $f_2(x,y)=0$ if $f_1(p)=f_2(p)=0$. An
intersection point must be always counted as times as its
multiplicity. The multiplicity of a point of intersection is a
rather complicated notion, which can be found, for instance, in
page $60$ of \cite{Walker}. Intuitively, the multiplicity of a
point of intersection takes into account the number of tangents
shared by the two curves at that point.
\par
B\'ezout's theorem takes into account the degree of two curves and
their points of intersection and relates them with the fact of
having a common factor.
\begin{theorem} {\sc [B\'ezout]} \ If $f_1(x,y)=0$ and $f_2(x,y)=0$ are two curves of degrees $m$ and
$n$, respectively, with more than $m n$ intersection points, then
there is an irreducible polynomial $r(x,y)$ which divides both
$f_1(x,y)$ and $f_2(x,y)$. \label{th4}
\end{theorem}

We are going to recall some results on fractionary power series
that can be found in \cite{Casas}. If $x$ is any free variable
over $\mathbb{C}$, we denote by $\mathbb{C}((x))$ the field of
fractions of $\mathbb{C}[[x]]$, where $\mathbb{C}[[x]]$ is the
ring of entire formal power series in $x$ with complex
coefficients. We recall that given a ring, which needs to be an
integral domain, we can define its field of fractions as the
smallest field containing it. The field of fractions is therefore
obtained from the integral domain by adding the least needed to
make of it a field, that is, the possibility of dividing by any
nonzero element. Given $n \in \mathbb{N}$, we consider entire
fractionary series $s= \sum_{i \geq r} a_i x^{i/n}$ where $r \in
\mathbb{Z}$, $a_i \in \mathbb{C}$ and $\min \{ i  :  a_i \neq 0\}
\geq 0$. An element $s \in \mathbb{C}((x^{1/n}))$ is of the form
$s= \sum_{i \geq r} a_i x^{i/n}$ where $r \in \mathbb{Z}$, $a_i
\in \mathbb{C}$ and the $\min \{ i  :  a_i \neq 0\}$ can be lower
than zero. The elements of the ring $\mathbb{C}[[x^{1/n}]]$ are
the entire fractionary series such that $\min \{ i  :  a_i \neq
0\} \geq 0$. It can be shown (see pages 17 and 18 of
\cite{Casas}), that given a fractionary series $s= \sum_{i \geq r}
a_i x^{i/n} \in \mathbb{C}((x^{1/n}))$ we can always take an
equivalent series such that $n$ and $\gcd \{ i : a_i \neq 0 \}$
have no common factor. In this case, we say that $n$ is the {\em
polydromy order} of the fractionary series $s$ and we denote it by
$\nu(s)$. Let $M(x,y)$ be a polynomial in $y$ of degree $N \in
\mathbb{N}$ whose coefficients in $y$ are fractionary series in
$\mathbb{C}((x^{1/n_i}))$, that is, we expand $M(x,y)$ in powers
of $y$: $M(x,y)=\sum_{i=0}^{N} s_i(x) y^i$ and we have that
$s_i(x) \in \mathbb{C}((x^{1/n_i}))$, for $i=0,1,2, \ldots, N$.
The polydromy order of $M(x,y)$ is defined as the least common
multiple of the polydromy orders of $s_i(x)$, $i=0,1,2,\ldots, N$,
that is, $\nu(M)=\mbox{\rm lcm} \{\nu(s_i) : i=0,1,2,\ldots, N\}$.
\par
The fractionary series $s$ is said to be {\em convergent} if
$\sum_{i \geq r} a_i t^{i}$ has non-zero convergence radius, where
$t=x^{1/n}$.
\par
The ring $\mathbb{C}\{x,y\} $ is the ring of convergent power
series in two variables and complex coefficients. The following
result clarifies the structure of an algebraic function and it is
stated and proved in \cite{Casas} (page $26$).
\begin{theorem}
If $f(x,y) \in \mathbb{C}[x,y]$, then there is a unique
decomposition of the form $ f = u x^r \prod_{i=1}^{\ell}
(y-g_i(x)), $ where $r \in \mathbb{N} \cup \{0 \}$, $u$ is an
invertible power series of  $\mathbb{C}\{ x,y \}$ and $g_i(x)$ are
convergent fractionary series. \label{th5}
\end{theorem}
In the book \cite{Casas} this result is stated for formal power
series $f(x,y) \in \mathbb{C}[[x,y]]$. In the Corollary 1.5.6 of
page $26$ in \cite{Casas}, it is stated that $\ell$ is the height
of the Newton polygon associated to $f(x,y)$. Since we are
considering a polynomial $f(x,y)$ instead of a formal series, we
deduce that $\ell$ is the highest degree in $y$ of $f(x,y)$ and,
therefore, $u$ is an invertible power series of $\mathbb{C} \{x
\}$. We define $h(x) = u x^r$ and we have that given $f(x,y) \in
\mathbb{C}[x,y]$ of degree $\ell$ in $y$, there is a unique
decomposition of the form:
\begin{equation}
f = h(x) \, \prod_{i=1}^{\ell} (y-g_i(x)), \label{eq4}
\end{equation}
where $h(x) \in \mathbb{C} \{ x \}$ and $g_i(x)$ are fractionary
series. Theorem 1.7.2 (page $31$) of \cite{Casas} ensures that all
the $y$--roots of $f(x,y)$ are convergent. Hence, we deduce that
an algebraic function $g(x)$ is a convergent fractionary series.
\par
By definition, given an algebraic function $g(x)$, we have that
there exists a polynomial such that $f(x, g(x)) \equiv 0$. If
$f(x,y)$ is not irreducible, then $g(x)$ is a fractionary series
appearing in the decomposition (\ref{eq4}) of at least one of the
irreducible factors of $f(x,y)$. Hence, without loss of
generality, we may always assume that $f(x,y)$ is irreducible.
Moreover, given $g(x)$, this irreducible polynomial $f(x,y)$ is
unique (modulus multiplication by constants). This statement is
clear from the fact that if $f_1(x,y)$ and $f_2(x,y)$ are two
irreducible polynomials such that $f_i(x,g(x)) \equiv 0$, $i=1,2$,
then these two polynomials have an infinite number of points of
intersection (because $g(x)$ is a convergent fractionary series)
and, by B\'ezout's Theorem \ref{th4}, we have that $f_1(x,y)$ and
$f_2(x,y)$ must be associates.

\section{The Main Results \label{sect3}}

\subsection{Particular algebraic solutions \label{sect31}}

\begin{theorem}
Let $g(x)$ be an algebraic particular solution of equation {\rm
(\ref{eq2})} and we call $f(x,y)$ the irreducible polynomial
satisfying $f(x,g(x)) \equiv 0$. Then, the curve $f(x,y)=0$ is an
invariant algebraic curve of system {\rm (\ref{eq1})}. \label{th6}
\end{theorem}
{\em Proof.} Let us denote by $F(x,y)$ the polynomial in
$\mathbb{R}[x,y]$ defined by: \[ F(x,y):= P(x,y)\frac{\partial
f}{\partial x} (x,y) + Q(x,y) \frac{\partial f}{\partial y}(x,y).
\] We have that $F(x,y)=0$ and $f(x,y)=0$ intersect in all the
points of the form $(x,g(x))$ by virtue of (\ref{eq3}). By
B\'ezout's Theorem \ref{th4}, we deduce that the polynomials
$f(x,y)$ and $F(x,y)$ share a common factor, because they
intersect in an infinite (continuum) number of points. They
intersect in all the points $(x,g(x))$ where the fractionary power
series $g(x)$ is convergent. Since $f(x,y)$ is an irreducible
polynomial, we have that $f(x,y)$ divides $F(x,y)$ in the ring of
real polynomials. From this fact, we conclude that there exists a
polynomial $k(x,y)$ such that $F(x,y)=k(x,y) f(x,y)$ and we get
that $f(x,y)=0$ is an invariant algebraic curve of system
(\ref{eq1}). \bbox \\

\noindent We recall the definition of invariant and of
quasipolynomial cofactor stated in \cite{GaGG}.
\begin{definition}
An {\em invariant} of {\rm (\ref{eq1})} is a function $\phi(x,y)$
such that there exists a {\rm quasipolynomial cofactor} $M(x,y)$,
where $M(x,y)$ is a polynomial in one of the variables $x$ or $y$
of degree $\leq m-1$ with $m$ the degree of system {\rm
(\ref{eq1})} in that variable, satisfying $P (\partial
\phi/\partial x)+ Q (\partial \phi/\partial y) = M \phi$.
\label{def7}
\end{definition}
In case that the set of points in $\mathbb{C}^2$ satisfying that
$\phi(x,y) = 0$ is not null, we have that $\phi$ is an invariant
if, and only if, $P (\partial \phi/\partial x)+ Q (\partial
\phi/\partial y)_{ | \, \phi=0 } =0$. We say that $\phi=0$ is an
{\em invariant curve} in this case.
\par
These definitions are a generalization of the so called
generalized cofactor introduced in \cite{GaG} where a
generalization of the Darboux integrability theory in order to
find non-Liouvillian first integrals of system (\ref{eq1}) was
presented. For the special invariant curve $\phi(x,y) := y - g(x)
= 0$ of (\ref{eq1}), where $g(x)$ is a particular solution of
equation (\ref{eq2}), a quasipolynomial cofactor always exists as
it was established in \cite{GaGG}. Many examples of invariants
with a quasipolynomial cofactor are given in \cite{GGG} as well as
a method to find first integrals, which are non-Liouvillian in
general, for certain families of systems.
\begin{proposition} {\sc \cite{GaGG}} \
A particular solution $g(x)$ of equation {\rm (\ref{eq2})} always
has a unique associated quasipolynomial cofactor of the form
$M(x,y) = k_{m-1}(x) y^{m-1} + \cdots + k_1(x) y + k_0(x)$, where
$m$ is the degree in $y$ of system {\rm (\ref{eq1})}.
\par
In case that $g(x)$ is algebraic then each $k_i(x)$, $i=0,1,2,
\ldots, m-1$ is a rational function in $x$ and $g(x)$ with
coefficients in $\mathbb{C}$. \label{prop8}
\end{proposition}
The second part of this proposition is deduced from the proof of
the first part as given in \cite{GaGG}.
\par
In case $g(x)$ is an algebraic function, by Theorem \ref{th5}, we
have that it is a fractionary series. Since $k_i(x)$ is a rational
function in $x$ and $g(x)$ with coefficients in $\mathbb{C}$ and
$g(x) \in \mathbb{C}((x^{1/n}))$, for a certain natural $n$, then
its quasipolynomial cofactor $M(x,y) \in
\mathbb{C}((x^{1/n}))[y]$, that is, each of the functions
$k_i(x)$, $i=0, 1,2,\ldots, m-1$ is a fractionary series with a
polydromy order that divides the polydromy order of $g(x)$. Next
lemma shows that in case $M(x,y)$ is a polynomial, then $g(x)$
must be a rational function.
\par
We say that an equation (\ref{eq2}) is {\em linear} if it is of
the form \[ \frac{dy}{dx} = m_1(x) \, y + m_0(x),
\] where $m_0(x)$ and $m_1(x)$ are functions of $x$. If equation
(\ref{eq2}) does not take this form, we say that it is {\em
non-linear}.
\begin{lemma}
Let $g(x)$ be a particular solution of a rational non-linear
equation {\rm (\ref{eq2})}. If the quasipolynomial cofactor
$M(x,y)$ related to $g(x)$ is a polynomial, then $g(x)$ is a
rational function. \label{lem9}
\end{lemma}
\par
{\em Proof.} We have that $M(x,y) = k_{m-1}(x) y^{m-1} + \cdots +
k_1(x) y + k_0(x)$ is a polynomial, so, $k_j(x)$ is a polynomial
in $x$ for all $j=0,1,2,\ldots, m-1$. Let us expand the
polynomials $P(x,y)$ and $Q(x,y)$ in powers of $y$: $P(x,y)=p_0(x)
+ p_1(x) y + \ldots + p_m(x) y^m$, $Q(x,y)=q_0(x) + q_1(x) y +
\ldots + q_m(x) y^m$, where $p_j(x)$ and $q_j(x)$ are the
polynomials in $x$ corresponding to the coefficients of degree $j$
in $y$ of $P(x,y)$ and $Q(x,y)$, respectively. Since $g(x)$ is a
particular solution of equation (\ref{eq2}), we have that $-P(x,y)
g'(x) + Q(x,y) = M(x,y)(y-g(x))$. Equating the coefficients of
order $j$ in $y$, we deduce that:
\[ p_0(x) g'(x) - k_0(x) g(x) = q_0(x), \quad (eq_0) \]
\[ p_j(x) g'(x) - k_j(x) g(x) = q_j(x) - k_{j-1}(x), \quad (eq_j) \ j=1,2,\ldots, m-1, \]
\[ p_m(x) g'(x) = q_m(x) - k_{m-1}(x) \quad (eq_m). \]
If $k_0(x)p_j(x)-p_0(x)k_j(x) \not\equiv 0$ for some
$j=1,2,\ldots, m-1$, we can equate $g(x)$ from the equations
$(eq_0)$ and $(eq_j)$ and we get \[ g(x) = \frac{p_0(x) q_j(x) -
q_0(x) p_j(x) -p_0(x) k_{j-1}(x)}{k_0(x)p_j(x)-p_0(x)k_j(x)}, \]
which is a rational function. If $k_0(x)p_j(x)-p_0(x)k_j(x) \equiv
0$ for all $j=1,2,\ldots, m-1$, we deduce that $p_j(x) = k_j(x)
L_1(x)$ and $p_0(x)=k_0(x) L_1(x)$ for all $j=1,2,\ldots, m-1$,
where $L_1(x)$ is a rational function in $x$.
\par
If $p_m(x) k_0(x) \not\equiv 0$, from the first and the last
equations we deduce that \[ g(x) = \frac{p_0(x) (q_m(x) -
k_{m-1}(x)) - p_m(x) q_0(x)}{p_m(x) k_0(x)}, \] which is a
rational function.
\par
We can also try to equate $g(x)$ from the equations $(eq_m)$ and
$(eq_j)$ for some $j=1,2,\ldots, m-1$. If $p_m(x) k_j(x)
\not\equiv 0$, then \[ g(x) = \frac{ q_m(x) p_j(x) - p_m(x) q_j(x)
+ p_m(x) k_{j-1}(x) - k_{m-1}(x) p_j(x)}{p_m(x) k_j(x)}, \] which
is a rational function.
\par
In case that $k_0(x)p_j(x)-p_0(x)k_j(x) \equiv 0$, $p_m(x) k_0(x)
\equiv 0$ and $p_m(x) k_j(x) \equiv 0$ for all $j=1,2,\ldots,
m-1$, first assume that $M(x,y) \equiv 0$. Then $y-g(x)$ would be
a first integral for system (\ref{eq1}), which means that $Q(x,y)
-P(x,y) g'(x) \equiv 0$. Hence, $Q(x,y)$ and $P(x,y)$ are such
that equation (\ref{eq2}) is a linear one, in contradiction with
our hypothesis. So, we conclude that $p_m(x) \equiv 0$. We also
have that $P(x,y)=L_1(x) M(x,y)$, that is, $p_j(x) = k_j(x)
L_1(x)$ for $j=0,1,2,\ldots, m-1$. From equation $(eq_m)$ we have
that $q_m(x)=k_{m-1}(x)$ and the other equations read for:
\[ k_0(x) (L_1(x) g'(x) - g(x)) = q_0(x), \quad (eq_0'), \]
\[ k_j(x) (L_1(x) g'(x) - g(x)) = q_j(x) - k_{j-1}(x), \quad
(eq_j'), \] for $j=1,2,\ldots, m-1$. Equating the factor $(L_1(x)
g'(x) - g(x))$ from equation $(eq_0')$ and $(eq_j')$ we deduce
that $k_j(x) q_0(x) - k_0(x) q_j(x) + k_0(x) k_{j-1}(x) \equiv 0$
for all $j=1,2,\ldots, m-1$. We write $q_0(x)=L_0(x) k_0(x)$,
where $L_0(x)$ is a rational function, and we have that $k_0(x)
(k_j(x) L_0(x) -q_j(x)+k_{j-1}(x)) \equiv 0$ for all
$j=1,2,\ldots, m-1$. If $k_0(x) \equiv 0$, then $P(x,y)$ and
$Q(x,y)$ share the common factor $y$. If $k_0(x) \not\equiv 0$,
then we have that $q_j(x) = k_j(x) L_0(x) + k_{j-1}(x)$, from
which we deduce that $Q(x,y)=(y+L_0(x)) M(x,y)$. \par Unless
$M(x,y)$ is a real number (different from zero), we have that
$P(x,y)$ and $Q(x,y)$ share a common factor, in contradiction with
our hypothesis. Therefore, after a rescaling of the time, the only
systems of the form (\ref{eq1}) with a particular non rational
solution $g(x)$ with a polynomial cofactor (in fact, the cofactor
is the real number $1$) are of the form : $\dot{x} = L_1(x)$,
$\dot{y}=y+L_0(x)$, which give rise to a linear equation. \bbox
\par
In the same way as in Lemma \ref{lem9}, we have that if $g(x)$ is
an algebraic particular solution of a nonlinear equation
(\ref{eq2}) with polydromy order $\nu(g)$, then its
quasipolynomial cofactor must have the same polydromy order.
\begin{lemma}
Let $g(x)$ be an algebraic particular solution of a non--linear
equation {\rm (\ref{eq2})}. Its quasipolynomial cofactor $M(x,y)$
has the same polydromy order as $g(x)$. \label{lem10}
\end{lemma}
{\em Proof.} We have that if $g(x)$ has polydromy order $\nu(g)$,
then the polydromy order of its quasipolynomial cofactor $M(x,y)$,
$\nu(M)$, divides it, that is, $\nu(M) | \nu(g)$. This is because
$M(x,y)$ is a polynomial in $y$ and a rational function in $x$ and
$g(x)$ as it has been stated in Proposition \ref{prop8}. \par The
fact that both polydromy orders coincide is a corollary of the
previous proof of Lemma \ref{lem9}. The reasonings are the same
since we can equate $g(x)$ in terms of the coefficients of
$M(x,y)$ unless the equation is of linear type. \bbox
\par
\begin{remark} If we have a linear equation of the form
$dy/dx \, = \, m_1(x) \, y + m_0(x)$, we may have a nonrational
particular solution with a polynomial cofactor. For instance,
taking $m_0(x) \equiv 0$ and $m_1(x) \equiv 1$, we have that
$g(x):=e^x$ is a particular solution with cofactor $1$. In the
same way, the linear equation $dy/dx \, = \, 3 y /(2 x)$ has the
algebraic solution $y= x^{3/2}$, whose polydromy order is $2$,
with the polynomial cofactor $3$, whose polydromy order is $1$.
\label{rem11}
\end{remark}

\subsection{On the invariants giving rise to exponential factors \label{sect32}}
In the Darboux theory of integrability not only invariant
algebraic curves are considered, but also exponential factors, as
we have stated in the introduction. Exponential factors appear in
a natural way from the coalescence of invariant algebraic curves,
as it is explained in \cite{Christopher1}. This statement means
that if we have a polynomial system (\ref{eq1}) with an
exponential factor of the form $\exp \{h / f_0 \}$, then there is
a $1$--parameter perturbation of system (\ref{eq1}), given by a
small $\varepsilon$, with two invariant algebraic curves, namely
$f_0=0$ and $f_0 + \varepsilon \, h=0$. Hence, when $\varepsilon
=0$, these two curves coalesce giving the exponential factor $\exp
\{h / f_0 \}$ for the system with $\varepsilon =0$, as well as the
invariant algebraic curve $f_0=0$ which does not disappear.
\par
In this context, and taking algebraic particular solutions into
account, the following question arises: which kind of function
appears with the coalescence of two algebraic particular
solutions? In the next Proposition \ref{prop12}, we show that the
natural generalization of algebraic particular solutions in this
framework is a function of the form: $\Phi(x,y):= \exp \{h_2(x) \,
A_1(x,y) / A_0(x,y) \}$ where $A_1(x,y)= \prod_{k=1}^{r}
(y-a_k(x))$ and $A_0(x,y)= \prod_{j=1}^{s} (y-\tilde{g}_j(x))$,
with $a_k(x)$ and $\tilde{g}_j(x)$ algebraic functions,
$k=1,2,\ldots,r$ and $j=1,2,\ldots, s$, $A_1(x,y)$ and $A_0(x,y)$
do not share any common factor, and $h_2(x)$ is a rational
function in $x$. For convention, if $r=0$ or $s=0$, we mean that
$A_1(x,y)$ or $A_0(x,y)$ takes a constant value, respectively,
which we may assume to be equal to $1$.
\par
Since by Proposition \ref{prop8}, we always have that a particular
solution given by $\phi(x,y)=y-g(x)$ has a quasipolynomial
cofactor, the property of being an invariant for $\Phi(x,y)$ is
given by associating to it a quasipolynomial cofactor $M(x,y)$
which is a polynomial in $y$ of degree at most $m-1$, where $m$ is
the degree in $y$ of system (\ref{eq1}). Moreover, we will assume
that $M(x,y)$ is well defined over $A_0(x,y)=0$, that is,
$M(x,\tilde{g}_j(x))$ is a real function of $x$ for all
$j=1,2,\ldots, s$. The following proposition gives the
characterization of this fact.
\begin{proposition}
We consider a function $\Phi(x,y):= \exp \{h_2(x) \, A_1(x,y) /
A_0(x,y) \}$ where $A_1(x,y)= \prod_{k=1}^{r} (y-a_k(x))$ and
$A_0(x,y)= \prod_{j=1}^{s} (y-\tilde{g}_j(x))$, with $a_k(x)$ and
$\tilde{g}_j(x)$ algebraic functions, $k=1,2,\ldots,r$ and
$j=1,2,\ldots, s$ $A_1(x,y)$ and $A_0(x,y)$ do not share any
common factor, and $h_2(x)$ is a rational function in $x$. Assume
that there exists a function $M(x,y)$, which is a polynomial in
$y$ of degree at most $m-1$, where $m$ is the degree in $y$ of
system {\rm (\ref{eq1})}, such that:
\begin{equation}
P(x,y) \frac{\partial \Phi(x,y)}{\partial x} + Q(x,y)
\frac{\partial \Phi(x,y)}{\partial y} = M(x,y) \, \Phi(x,y) .
\label{eq5}
\end{equation}
We assume that $M(x,\tilde{g}_j(x))$ is a real function of $x$ for
all $j=1,2,\ldots, s$. We denote by $\tilde{k}_i(x)$ the
coefficient of degree $i$ in $y$ of $M(x,y)$, that is, $M(x,y) =
\tilde{k}_0(x) + \tilde{k}_1(x) y + \tilde{k}_2(x) y^2 + \ldots +
\tilde{k}_{m-1}(x) y^{m-1}$. Then,
\begin{itemize}
\item[{\rm (i)}] Each one of the algebraic functions
$\tilde{g}_j(x)$, $j=1,2,\ldots, s$ is a particular solution of
equation {\rm (\ref{eq2})}. We denote by $M_j(x,y)$ its associated
quasipolynomial cofactor.

\item[{\rm (ii)}] The following identity is satisfied:
\begin{equation}
\mathcal{X} \Big(h_2(x)\, A_1(x,y) \Big) = \left( \sum_{j=1}^{s}
M_j(x,y) \right) h_2(x)\, A_1(x,y) \, + \, M(x,y) \, A_0(x,y) .
\label{eq6}
\end{equation}

\item[{\rm (iii)}] Each one of the functions $\tilde{k}_i(x)$ is
rational in $x$ and rational in $a_k(x)$ and $\tilde{g}_j(x)$,
with $k=1,2,\ldots,r$ and $j=1,2,\ldots, s$.

\end{itemize}
\label{prop12}
\end{proposition}
{\em Proof.} We have that $\mathcal{X} \left( \Phi(x,y) \right) =
M(x,y) \Phi(x,y)$, from which we deduce the following identity:
\begin{equation} \mathcal{X} \left( h_2(x)\, A_1(x,y) \right) \, A_0(x,y) -
h_2(x) A_1(x,y) \, \mathcal{X} \left( A_0(x,y) \right) = M(x,y) \,
A_0(x,y)^2 . \label{eq7}
\end{equation}
Since $A_1(x,y)$ and $A_0(x,y)$ do not share any common factor and
$A_0(x, \tilde{g}_j(x)) \equiv 0$ for $j=1,2,\ldots, s$, then,
from equation (\ref{eq7}), we have that $\mathcal{X}\left(
A_0(x,y) \right)|_{y=\tilde{g}_j(x)} \equiv 0$. We notice that
here we are using that $M(x,\tilde{g}_j(x))$ is a real function of
$x$ for all $j=1,2,\ldots, s$. Let us call $A_{0_j}(x,y)
:=\prod_{i=1, i\neq j}^{s} (y-g_i(x))$ and we have that
$A_0(x,y)=A_{0_j}(x,y) (y-\tilde{g}_j(x))$. Then,
\[ \mathcal{X}\left( A_0(x,y) \right) = \mathcal{X}\left( A_{0_j}
(x,y) \right) (y-\tilde{g}_j(x)) + A_{0_j}(x,y) \left(
-\tilde{g}_j'(x) P(x,y) + Q(x,y) \right).
\] $ $From $\mathcal{X}\left( A_0(x,y)
\right)|_{y=\tilde{g}_j(x)} \equiv 0$, we deduce that
$-\tilde{g}_j'(x) P(x,\tilde{g}_j(x)) + Q(x,\tilde{g}_j(x)) \equiv
0$ and we conclude that $\tilde{g}_j(x)$ is a particular solution
of system (\ref{eq1}), as stated in (i).
\par
In order to prove (ii), we consider the quasipolynomial cofactor
$M_j(x,y)$ associated to $y-\tilde{g}_j(x)$, whose existence is
ensured by Proposition \ref{prop8}. Then, we have that
$\mathcal{X} (A_0) = \left( \sum_{j=1}^{s} M_j \right) \, A_0$.
Hence, identity (\ref{eq7}) reads for: \[ \begin{array}{l}
\displaystyle \mathcal{X} \Big( h_2(x)\, A_1(x,y) \Big) \,
A_0(x,y) - \, h_2(x) \, A_1(x,y) \, \left( \sum_{j=1}^{s} M_j(x,y)
\right) A_0(x,y) \ =
\\ \displaystyle \qquad =  M(x,y) \, A_0(x,y)^2 .
\end{array} \] This identity coincides with (\ref{eq6}), after
dividing both members by $A_0(x,y)$.
\par Finally, to prove (iii), we observe that equating the
coefficients of the same degree in $y$ from identity (\ref{eq6}),
we deduce that each $\tilde{k}_i(x)$ is a rational function of
$x$, $\tilde{g}_j(x)$, $a_k(x)$ and $a_k'(x)$. Since each $a_k(x)$
is an algebraic function, we can consider $f_{a_k}(x,y)$ as the
irreducible polynomial in $\mathbb{C}[x,y]$ such that
$f_{a_k}(x,a_k(x)) \equiv 0$. We can derive this last identity
with respect to $x$ and we deduce that: \[ \frac{\partial
f_{a_k}(x,a_k(x))}{\partial x}\, + \, \frac{\partial
f_{a_k}(x,a_k(x))}{\partial y} \  a_k'(x) \ \equiv 0. \] Hence, we
can substitute the value $a_k'(x)$ appearing in $\tilde{k}_i(x)$
by the rational expression in $x$ and $a_k(x)$: $- \left[
\partial f_{a_k}(x,a_k(x))/\partial x \right] / \left[
\partial f_{a_k}(x,a_k(x))/\partial y \right]$. Therefore, each
$\tilde{k}_i(x)$ is a rational function of $x$, $\tilde{g}_j(x)$
and $a_k(x)$, as we wanted to prove. \bbox
\newline

As Theorem \ref{th6} shows, an algebraic particular solution
$y-g(x)$ recovers the invariant algebraic curve to which it is
related. That is, we have that $f(x,y)=0$ is an invariant
algebraic curve of system (\ref{eq1}) if, and only if, all its
$y$--roots are algebraic particular solutions of equation
(\ref{eq2}). The fact of being particular solutions implies that
each one of them has an associated quasipolynomial cofactor. We
notice that all the $y$--roots appear in the factorization
$f(x,y)=h(x) \prod_{i=1}^{\ell} (y-g_i(x))$ described by
(\ref{eq5}). We would like to have an analogous to this statement
but for exponential factors and this is what is given in the next
Theorem \ref{th13}.
\par
In this context, if we have an exponential factor given by
$\exp\{R(x,y)\}$, where $R(x,y)$ is a rational function, i.e.
$R(x,y) \in \mathbb{R}(x,y)$, we would like to write a
factorization for it analogous to (\ref{eq5}). We notice that
given any two functions $R_1(x,y)$, $R_2(x,y)$ with the property
that $R_1(x,y) + R_2(x,y) \equiv R(x,y)$, we have that
$\exp\{R(x,y)\} \equiv \exp \{R_1(x,y)\} \cdot \exp \{R_2(x,y)
\}$, which is a ``factorization'' of $\exp\{R(x,y)\}$. However,
not all these factorizations are useful to our purposes because it
can be shown that $\exp \{R_1(x,y)\}$ and $\exp \{R_2(x,y)\}$ do
not need to have an associated quasipolynomial cofactor. So, they
are not invariants for system (\ref{eq1}), as the exponential
factor $\exp\{R(x,y)\}$ is. The following example exhibits this
fact.
\par
{\bf Example.} Let us consider the following cubic system with the
invariant straight line $y=0$ :
\begin{equation}
\dot{x} \, = \, (2 x+y)(1+x) + 2 x^2 y + y^3, \qquad \dot{y}= y
(1+x+x y). \label{eqex}
\end{equation}
This system has the invariant $\Phi(x,y):= \exp \{ \sqrt{x} / y
\}$ with the quasipolynomial cofactor $M(x,y)= (1+x+y^2)/(2
\sqrt{x})$. As it will be proved in Theorem \ref{th13}, the
existence of this invariant implies the existence of the following
exponential factor: $F(x,y):=\exp \{ (x+y)/y^2 \}$ with the
cofactor $\tilde{k}_0 (x,y)=y-x$. We notice that not all the
functions $R_1, R_2$ satisfying $R_1+R_2= (x+y)/y^2$ give rise to
an invariant. For instance, if we take $R_1:=x/y^2$ and
$R_2:=1/y$, the following easy computation shows that $\exp
\{R_1\}$ is not an invariant of system (\ref{eqex}), and neither
$\exp \{R_2\}$ is. We note that, taking $\mathcal{X}$ as the
vector field associated to system (\ref{eqex}), we have that
$\mathcal{X} (\exp\{R_1\}) = (y + (1+x)/y) \, \exp\{R_1\}$ and
since the rational function $(y + (1+x)/y)$ is not well-defined
over $y=0$, we deduce that $\exp\{R_1\}$ is not an invariant.
However, if we consider $\Psi(x,y):= \exp \{ (x+y-y \sqrt{x})/y^2
\}$, we have that $\Phi(x,y) \cdot \Psi(x,y) = F(x,y)$ and both
$\Phi(x,y)$ and $\Psi(x,y)$ are invariants of system (\ref{eqex}).
\newline

The following theorem shows that if $\Phi(x,y):= \exp \{h_2(x)\,
A_1(x,y) / A_0(x,y) \}$, as described in Proposition \ref{prop12},
is such that it has an associated quasipolynomial cofactor, then
there exists another function $\Psi(x,y):=\exp \{B_1(x,y) /
B_0(x,y)\}$, of the same type as $\Phi(x,y)$, such that the
product $\Phi(x,y) \Psi(x,y)$ gives rise to an exponential factor
$\exp\{R(x,y)\}$ for system (\ref{eq1}). This sentence means that
a factorization of $\exp\{R(x,y)\}$ giving invariants for system
(\ref{eq1}) is obtained by the product $\Phi(x,y) \Psi(x,y)$.
Therefore, we are able to recover $\exp\{R(x,y)\}$ from one of its
factors $\Phi(x,y)$ and, moreover, since $\exp\{R(x,y)\}$ appears
by the coalescence of two invariant algebraic curves, we conclude
that $\Phi(x,y)$ needs to appear by the coalescence of algebraic
particular solutions since it is formed by a product of them. This
fact is the result that we were targeting to: we wanted to exhibit
the analogy between the generalization of invariant algebraic
curves to exponential factors with the generalization of algebraic
particular solutions to invariants of the form $\Phi(x,y):= \exp
\{h_2(x)\, A_1(x,y) / A_0(x,y) \}$.
\begin{theorem}
Assume that the function $\Phi(x,y):= \exp \{h_2(x)\, A_1(x,y) /
A_0(x,y) \}$ has a quasipolynomial cofactor $M(x,y)$, that is,
$\mathcal{X} (\Phi(x,y)) = M(x,y) \, \Phi(x,y)$ and assume that
$M(x,y)$ is  well defined over $\{(x,y) \, : \, A_0(x,y)=0 \}$.
Then, there exist quasipolynomial functions $B_0(x,y)$ and
$B_1(x,y)$, which are polynomials in $y$ and algebraic in $x$,
such that $R:= h_2 \, A_1 /A_0 \, + \, B_1 / B_0$ is a rational
function in $x$ and $y$ and $\exp \{ R(x,y) \}$ is an exponential
factor of system {\rm (\ref{eq1})}. \label{th13}
\end{theorem}
{\em Proof.} We consider $\phi_A (x,y, \varepsilon):= A_0(x,y) +
\varepsilon h_2(x) \, A_1(x,y)$ which is a polynomial in $y$ whose
coefficients are algebraic functions. We can compute the sequence
of powers $\phi_A(x,y,\varepsilon)^{j+1}$, for each natural number
$j$ and the coefficients of all these polynomials in $y$ are
algebraic functions, combination of those of
$\phi_A(x,y,\varepsilon)$. Therefore, there exists a natural
number $N$ such that $\phi_A(x,y,\varepsilon)^{N+1}$ is a linear
combination of all the previous powers and the coefficients of
this combination are powers of $x$, $y$ and $\varepsilon$. This
fact is due to the finiteness of the algebraic extensions given by
the $y$--roots of $\phi_A(x,y,\varepsilon)$. In this way, we have
that there exists an irreducible polynomial
$\mathcal{P}(x,y,\varepsilon)$ in $x$, $y$ and $\varepsilon$, such
that each of the $y$--roots of $\phi_A(x,y,\varepsilon)$ is an
$y$--root of $\mathcal{P}(x,y,\varepsilon)$. This polynomial
$\mathcal{P}(x,y,\varepsilon)$ is the minimal polynomial of the
$y$-roots of $\phi_A(x,y,\varepsilon)$, and it can be always taken
irreducible. We expand $\mathcal{P}(x,y,\varepsilon)$ in powers of
$\varepsilon$ and we denote by $R_i(x,y)$ its coefficient of
degree $i$ in $\varepsilon$, which is a polynomial in $x$ and $y$.
Let us consider the quotient $\mathcal{P}(x,y, \varepsilon) /
\phi_A(x,y,\varepsilon)$ which is a polynomial in $y$, denoted by
$\phi_B(x,y,\varepsilon)$. We expand $\phi_B(x,y,\varepsilon)$ in
powers of $\varepsilon$ and we denote by $B_i(x,y)$ its
coefficient of degree $i$ in $\varepsilon$, which is a polynomial
in $y$. Since $\mathcal{P} = \phi_A \phi_B$, we deduce that $R_0 =
A_0 B_0$ and $R_1 = A_0 B_1 + h_2 A_1 B_0$, equating the
coefficients of $\varepsilon^{0}$ and $\varepsilon^1$. Thus, $R:=
h_2 A_1 /A_0 \, + \, B_1 / B_0 = R_1/R_0$ is a rational function
in $x$ and $y$. In case $R_1 \equiv 0$ all the reasoning works
just taking the first $R_i \not\equiv 0$ with $i>0$.
\par
We only need to see that $\exp \{ R(x,y) \}$ is an exponential
factor of system {\rm (\ref{eq1})}. We have, from Proposition
\ref{prop12}, that $A_0$ is the product of algebraic particular
solutions $y-\tilde{g}_j(x)$. Hence $\mathcal{X}\left( A_0(x,y)
\right) = M_0(x,y) A_0(x,y)$, where $M_0(x,y)$ is the
quasipolynomial cofactor $\sum_{j=1}^{s} M_j(x,y)$, which is a
polynomial in $y$ of degree at most $m-1$, $m$ being the degree in
$y$ of system (\ref{eq1}). By the irreducibility of
$\mathcal{P}(x,y, \epsilon)$, we deduce that $R_0(x,y)$ is a power
of the lowest degree polynomial containing as $y$--roots all the
$\tilde{g}_j(x)$. Hence, by Theorem \ref{th6}, we have that
$R_0(x,y)=0$ is an invariant algebraic curve of system
(\ref{eq1}), that is, $\mathcal{X}(R_0) = k_0 R_0$, where
$k_0(x,y)$ is a polynomial in $x$ and $y$. Moreover, from $R_0=A_0
B_0$ and $\mathcal{X}\left( A_0  \right) = M_0 A_0$ we deduce that
$\mathcal{X}\left( B_0 \right) = (k_0 - M_0) B_0$. By Lemma
\ref{lem1}, we only need to show that there is a polynomial
$\tilde{k}_0(x,y)$ such that $\mathcal{X} \left( R_1 \right) = k_0
R_1 + \tilde{k}_0 R_0$. We recall that from Proposition
\ref{prop12} we have that $\mathcal{X} \left( h_2 A_1 \right) =
M_0 h_2 A_1 + M \, A_0$. We consider the polynomial $G(x,y):=
\mathcal{X}\left( R_1(x,y) \right) - k_0(x,y) R_1(x,y)$. We have:
\begin{eqnarray*}
G & = & \mathcal{X}\left( R_1 \right) - k_0 R_1
\\ & = & \mathcal{X}\left(A_0 \right) B_1 + A_0 \mathcal{X}\left( B_1
\right) + \mathcal{X}\left( h_2 A_1 \right) B_0 + h_2 A_1
\mathcal{X}\left(B_0 \right) - k_0 \left( A_0 B_1 + h_2 A_1 B_0
\right) \\ & = & M_0 A_0 B_1 + A_0 \mathcal{X}\left( B_1 \right) +
(M_0 h_2 A_1 + M  A_0)  B_0 + h_2 A_1 (k_0 - M_0) B_0
\\ & &
- \, k_0 (A_0 B_1 + h_2 A_1 B_0) \\ & = & \left[ M_0 B_1 +
\mathcal{X}\left( B_1 \right) + M \, B_0 - k_0 B_1 \right] \, A_0.
\end{eqnarray*}
We deduce that $G(x,\tilde{g}_j(x)) \equiv 0$ for all the
$\tilde{g}_j(x)$ appearing in $A_0$. We note that here we are
using the hypothesis that $M(x,y)$ is  well defined over $\{(x,y)
\, : \, A_0(x,y)=0 \}$. Since $R_0(x,y)$ is the lowest degree
polynomial with this property and $G(x,y)$ is a polynomial, we
have, by B\'ezout's Theorem, that $R_0(x,y)$ is a divisor of
$G(x,y)$ in the ring of polynomials $\mathbb{C}[x,y]$. So, there
exists a polynomial $\tilde{k}_0(x,y)$ such that $G(x,y) =
\tilde{k}_0(x,y) R_0(x,y)$ and, thus, $\mathcal{X} \left( R_1
\right) = k_0 R_1 + \tilde{k}_0 R_0$. Therefore, $\exp\{ R_1/R_0
\}$ is an exponential factor of system (\ref{eq1}). \bbox

\subsection{On the structure of the quasipolynomial cofactor \label{sect33}}

The following proposition gives the form of the quasipolynomial
cofactor associated to an invariant of the form $I(x,y)=h(x) \,
\prod_{i=1}^{\ell} (y-g_i(x))^{\alpha_i}$.
\par
We define as $m$ the degree of system (\ref{eq1}) in the variable
$y$ and we expand the polynomials $P(x,y)$ and $Q(x,y)$ in this
variable:
\[ P(x,y)=\sum_{i=1}^{m} p_i(x) y^i, \quad Q(x,y)=\sum_{i=1}^{m} q_i(x)
y^i .\]
\begin{proposition}
Let $I(x,y)=h(x) \prod_{i=1}^{\ell} (y-g_i(x))^{\alpha_i} $ be an
invariant of system {\rm (\ref{eq1})} with an associated
quasipolynomial cofactor $k(x,y):= k_0(x) + k_1(x) y + \ldots +
k_{m-1}(x) y^{m-1}$. Then, $p_m(x) h'(x)  \equiv 0$ and
\begin{equation}
\displaystyle k_j(x) =  p_j(x) \, \displaystyle
 \frac{h'(x)}{h(x)} \, \, +  \displaystyle
\, \displaystyle \sum_{s=j+1}^{m} \left( \sigma_{s-(j+1)}(x) \,
q_s(x) - \displaystyle \frac{\sigma_{s-j}'(x)}{(s-j)} \, p_s(x)
\right), \label{eq10}
\end{equation}
for $j=0,1,2, \ldots, m-1$, where
\[ \sigma_{\kappa}(x) \ = \ \sum_{\nu=1}^{\ell} \alpha_{\nu} \,
g_{\nu}^{\kappa}(x), \qquad \mbox{for }\kappa=0,1,2, \ldots, m. \]
We notice that the first term in {\rm (\ref{eq10})} does not
appear if $h(x)$ is a constant or the last term does not appear if
$p_m(x)$ is zero. \label{prop14}
\end{proposition}
{\em Proof.}  As it has been shown in \cite{GaGG}, each of the
factors $y-g_i(x)$ involved in the expression $I(x,y)$ is a
particular solution of system (\ref{eq1}). Moreover, a particular
solution $y-g(x)$ has a related quasipolynomial cofactor of degree
at most $m-1$ in $y$, therefore we have that $k(x,y)$ is a
polynomial in $y$ of degree at most $m-1$. We will deduce each
expression of $k_i(x)$ from the identity:
\begin{equation}
P(x,y) \frac{\partial I}{\partial x}(x,y) + Q(x,y) \frac{\partial
I}{\partial y}(x,y) = k(x,y) I(x,y). \label{eq8}
\end{equation}
We compute the partial derivatives of $I(x,y)$ in (\ref{eq8}) and
we divide each member of the resulting expression by $I(x,y)$. We
obtain:
\[ P(x,y) \left[ \frac{h'(x)}{h(x)} - \sum_{\nu=1}^{\ell}
\frac{\alpha_\nu g_{\nu}'(x)}{y-g_{\nu}(x)} \right] + Q(x,y)
\left[ \sum_{\nu=1}^{\ell} \frac{\alpha_\nu}{y-g_\nu(x)} \right] =
k(x,y). \] Then, we deduce:
\begin{equation}
\begin{array}{l}
\displaystyle P(x,y) \left[ \frac{h'(x)}{h(x)}
\prod_{\nu=1}^{\ell} (y-g_{\nu}(x)) - \sum_{\nu=1}^{\ell}
\alpha_\nu g_\nu'(x) \prod_{\mu=1, \mu \neq \nu}^{\ell}
(y-g_{\mu}(x)) \right] +
\\ \displaystyle + \, Q(x,y) \sum_{\nu=1}^{\ell} \alpha_\nu \left( \prod_{\mu=1, \mu
\neq \nu}^{\ell} (y-g_{\mu}(x)) \right) \, = \, k(x,y)
\prod_{\nu=1}^{\ell} (y-g_{\nu}(x)).
\end{array}
\label{eq9}
\end{equation}
The expression (\ref{eq9}) is an identity of polynomials in $y$ of
degree $m+\ell$. Let us consider the equality of coefficients of
degree $m+\ell$ in $y$: $p_{m}(x) h'(x)/h(x) =0$, which implies
that either $p_{m}(x) \equiv 0$ or $h(x)$ is a constant function.
\par
We define $eq_i$ as the equation resulting from identifying the
coefficients of $y^i$ in both members of equation (\ref{eq9}),
$i=0,1,2,\ldots, m + \ell -1$. From $eq_{\ell+m-1}$, we can equate
the expression of $k_{m-1}(x)$. Once we know this expression, from
$eq_{\ell + m -2}$, we can equate the expression of $k_{m-2}(x)$.
Once we know this function, from $eq_{\ell + m -3}$ we can equate
the expression of $k_{m-3}(x)$, and so on. Hence, in a recursive
way, from $eq_{\ell + j}$, we equate $k_{j}(x)$, where
$j=m-1,m-2,m-3, \ldots, 2,1,0$. It can be shown by induction that
these expressions are given by (\ref{eq10}). We notice that
\[ \frac{\sigma_{\kappa}'(x)}{\kappa} \ = \ \sum_{\nu=1}^{\ell}
\alpha_{\nu} g_\nu^{\kappa-1}(x) g_{\nu}'(x) , \] for
$\kappa=0,1,2,\ldots,m$.
\par
We substitute the values of $k_j(x)$ given in (\ref{eq10}) in
equation (\ref{eq9}) and we deduce that a polynomial of degree at
most $\ell -1$ in $y$ must be zero. We denote by $Pol(y)$ this
polynomial in $y$, which is also a function of $x$, $h(x)$ and
$g_1(x), g_2(x), \ldots, g_\ell(x)$. We fix an index $i$ such that
$1 \leq i \leq \ell$ and when we substitute $y$ by $g_i(x)$ in
$Pol(y)$, we get that:
\[ Pol(g_i(x)) \, = \, \prod_{s=1,s\neq i}^{\ell} (g_i(x) -
g_s(x)) \, \Big[ Q(x,g_i(x)) - P(x,g_i(x)) g_i'(x)  \Big] . \]
Since each $g_i(x)$ is a particular solution of system
(\ref{eq1}), we deduce that each $g_i(x)$ is a different root of
the polynomial $Pol(y)$. Then, we have  $\ell$ different roots of
a polynomial of degree at most $\ell-1$, so $Pol(y)$ is the null
polynomial. We conclude that once we have substituted the cofactor
as defined by (\ref{eq10}) in equation (\ref{eq9}) we have that
this equation is satisfied. \bbox
\newline

In the same way as in Proposition \ref{prop14}, we are going to
give the form of the quasipolynomial cofactor associated to an
invariant of the form \begin{equation} I(x,y) = \exp \left\{
h_2(x) \frac{A_1(x,y)}{A_0(x,y)} \right\} h_1(x)
\prod_{i=1}^{\ell} (y-g_i(x))^{\alpha_i}, \label{eqfi}
\end{equation} where we have defined $A_1(x,y):= \prod_{k=1}^{r}
(y-a_k(x))$ and $A_0(x,y):= \prod_{j=1}^{s} (y-\tilde{g}_j(x))$.
\par
The first thing that we notice is that, without loss of
generality, we can assume that $r=s$. In case that $r \neq s$, we
can consider the change of variable $y = 1/z$. After this change
(and a reparameterization of the time variable if necessary), we
obtain another polynomial system (\ref{eq1}) and the function $I$
is transformed to another one with the same structure but with
$r=s$. This last assertion is clear from the following equality:
\[ \begin{array}{lll} \displaystyle h_2(x) \frac{A_1(x,1/z)}{A_0(x,1/z)} & = & \displaystyle h_2(x)
\frac{\prod_{k=1}^{r} \left( \frac{1}{z}-a_k(x)
\right)}{\prod_{j=1}^{s} \left( \frac{1}{z}-\tilde{g}_j(x)
\right)} = \\ &= & \displaystyle  h_2(x) \frac{ \frac{1}{z^r}
\prod_{k=1}^{r} \left( 1-a_k(x)z \right)}{\frac{1}{z^s}
\prod_{j=1}^{s} \left( 1-\tilde{g}_j(x)z \right)} = \tilde{h}_2(x)
\frac{z^{s-r} \prod_{k=1}^{r} \left( z - \frac{1}{a_k(x)} \right)
}{\prod_{j=1}^{s} \left( z - \frac{1}{\tilde{g}_j(x)} \right)},
\end{array} \] where $\tilde{h}_2(x):= h_2(x) (-1)^{s-r} \prod_{k=1}^{r} a_k(x)
\diagup \prod_{j=1}^{s} \tilde{g}_j(x)$. If we ensure the
structure of a quasipolynomial cofactor for the system with
variables $(x,z)$, we deduce its structure for the system with
variables $(x,y)$ just undoing the change and the
reparameterization of the time, if it has been done. So, without
loss of generality, we can assume that $r=s$. In fact, an
analogous proof can be done for the case $r \neq s$ but as it
involves many computations which are not far from the ones that we
are showing, we have preferred to avoid the case $r \neq s$ by
means of the change of variable $y=1/z$.
\begin{proposition}
Let $I(x,y)$ be an invariant of system {\rm (\ref{eq1})} of the
form {\rm (\ref{eqfi})} with $r=s$ and with an associated
quasipolynomial cofactor $k(x,y):= k_0(x) + k_1(x) y + \ldots +
k_{m-1}(x) y^{m-1}$. Then, $p_m(x) [ h_1'(x) + h_1(x) h_2'(x) ]
\equiv 0$ and
\begin{equation}
\displaystyle k_j(x)  =  \displaystyle \left(
\frac{h_1'(x)}{h_1(x)} + h_2'(x) \right) p_j(x) \, + \,
\displaystyle \displaystyle \sum_{i=j+1}^{m} \left(
\tilde{\sigma}_{i-j-1}(x) \, q_{i}(x) -
\frac{\tilde{\sigma}_{i-j}'(x)}{(i-j)} \, p_{i}(x) \right),
\label{eq11a}
\end{equation}
for $j=0,1,2,\ldots, m-1$, where either the first term does not
appear in case that $h_1(x) = c\exp\{ - h_2(x) \}$, with $c \in
\mathbb{C}$ and $c \neq 0$, or the last term does not appear if
$p_m(x) \equiv 0$. The functions $\tilde{\sigma}_{\kappa}(x)$ are
defined in the following way. Given $\kappa \in \mathbb{N}$ we
consider the set of indexes
\[ J_{\kappa} := \left\{ (\epsilon_{1},\epsilon_{2},\ldots,
\epsilon_{r},i_1,i_2, \ldots, i_r) \, : \, \sum_{k=1}^{r}
\epsilon_{k} + \sum_{j=1}^{r} i_j = \kappa , \, \epsilon_j \in \{
0,1 \} , \, i_j \in \mathbb{N} \cup \{ 0 \} \right\} \] and we
have that:
\[ \tilde{\sigma}_{\kappa} := \sum_{\nu=1}^{\ell} \alpha_{\nu}
g_{\nu}^{\kappa} \, + \, \kappa \, h_2 \, \sum_{J_{\kappa}}
(-1)^{\epsilon_{1} + \epsilon_{2} + \ldots + \epsilon_{r} + 1} \,
a_1^{\epsilon_{1}} a_2^{\epsilon_{2}} \cdots a_r^{\epsilon_{r}}
\tilde{g}_1^{i_1} \tilde{g}_2^{i_2} \cdots \tilde{g}_s^{i_r}  ,\]
for $\kappa=0,1,2,\ldots, m$. \label{prop15}
\end{proposition}
{\em Proof.} We have that if $I(x,y)$ is an invariant, then each
one of the functions $g_i(x)$ and $\tilde{g}_j(x)$ are particular
solutions and the identity (\ref{eq6}) is satisfied. We are going
to deduce the expressions of $k_j(x)$ only assuming that the
following identity is satisfied:
\begin{equation} P(x,y) \frac{\partial I}{\partial x}(x,y) + Q(x,y) \frac{\partial
I}{\partial y}(x,y) = k(x,y) I(x,y), \label{eq11} \end{equation}
We multiply (\ref{eq11}) by \[ A_0(x,y)^2 \, \exp \left\{ -
\frac{h_2(x) A_1(x,y)}{A_0(x,y)} \right\} \, \prod_{i=1}^{\ell}
(y-g_i(x))^{1-\alpha_i} \] so as to get an identity of polynomials
in $y$ of degree $m + \ell + 2 r$, where $m$ is the degree of
system (\ref{eq1}) in $y$. The equality of the coefficients of
highest degree in $y$, that is $y^{m+ \ell + 2 r}$, gives $p_m(x)
[ h_1'(x) + h_1(x)  h_2'(x) ] \equiv 0$. This condition gives two
possibilities: either $p_m(x) \equiv 0$ or $h_1(x) = c \exp\{ -
h_2(x) \}$, with $c \in \mathbb{C}$ and $c \neq 0$. Since the
proof for both cases is analogous, we are going to follow them
simultaneously. It has been shown in \cite{GG} that $k(x,y)$ is a
polynomial in $y$ of degree at most $m-1$ so we can collect it in
this variable: $k(x,y)=\sum_{i=1}^{m-1} k_i(x) y^i$. The equality
of coefficients of degree $j+\ell + 2 r$ gives us the expression
of $k_j(x)$ in a recursive way. We first compute $k_{m-1}$ from
the equality of coefficients of $y^{m-1+\ell + 2 r}$, once we have
this one, we compute $k_{m-2}$ from the equality of coefficients
of $y^{m-2+\ell + 2 r}$ and so on. Some tedious computations show
that these expressions are the ones given in (\ref{eq11a}).
\par
We substitute the given expressions of $k_j(x)$ in the equality
(\ref{eq11}) and from the equality of the coefficients of $y^{\ell
+ 2 r-1}$ we can compute the function $h_2(x)$. Now we have a
polynomial in $y$ of degree at most $\ell+ 2 r -1$ in $y$ that
must be identically zero. We denote by $\overline{Pol}(y)$ this
polynomial in $y$. In the same way as in the proof of Proposition
\ref{prop14}, we substitute the $y$ variable by each one of the
$\ell + 2 r$ functions $g_i(x)$, $\tilde{g}_j(x)$ and $a_k(x)$ and
we deduce that $\overline{Pol}(g_i(x)) \equiv 0$ because $g_i(x)$
is a particular solution for $i=1,2,\ldots, \ell$,
$\overline{Pol}(\tilde{g}_j(x)) \equiv 0$ because $\tilde{g}_j(x)$
is a particular solution for $j=1,2,\ldots, r$ and
$\overline{Pol}(a_k(x)) \equiv 0$ because relation (\ref{eq6}) is
satisfied for $k=1,2, \ldots, r$. Therefore, we have $\ell+ 2 r$
different roots of a polynomial of degree at most $\ell+2r-1$, so
$\overline{Pol}(y)$ is the null polynomial. We conclude that once
we have substituted the cofactor defined by (\ref{eq11a}) and the
function $h_2(x)$ in equation (\ref{eq11}) we have that this
equation is satisfied. \bbox \par {\bf Example.} We include an
example of the form of the quasipolynomial cofactor of a function
$I(x,y)$ as given by Proposition \ref{prop15} deduced from the
equation (\ref{eq11a}). We consider the following system:
\begin{equation}
\dot{x}=y+y^2+x^2+4 y x^2, \quad \dot{y}=-x-2x^3+2 x y^2,
\label{eqex11}
\end{equation}
which has the exponential factor $I(x,y)=\exp\{ (2
y-1)/(x^2+y^2)\}$ with cofactor $k(x,y):=-4x$. We are going to
deduce the value of this cofactor by using the formulas stated in
(\ref{eq11a}) of Proposition \ref{prop15}. We notice that this
function $I(x,y)$ has $r=1<s=2$, so we need to perform the change
of variables $y=1/z$. The transformed system needs a
reparameterization of the time, meaning multiplying it by $z^2$,
in order to have a polynomial system. Then, the transformed system
reads for:
\begin{equation} \dot{x}=1+(1+4x^2)z + x^2 z^2, \quad \dot{z}=-2x z^2 +
(x+ 2 x^3) z^4, \label{eqex11b} \end{equation} which has the
exponential factor $I(x,z)=\exp\{ z(2-z)/(x^2 z^2+1)\}$ with
cofactor $k(x,z)=-4x z^2$. This cofactor coincides with the
transformation of the cofactor in coordinates $(x,y)$ after the
reparameterization of time. Using the notation described in
Proposition \ref{prop15} we have that, for system (\ref{eqex11b}),
$m=4$, $r=s=2$, $\ell=0$, $a_1(x):=0$, $a_2(x):=2$,
$h_2(x):=-1/x^2$, $\tilde{g}_1(x):=-i/x$, $\tilde{g}_2(x):=i/x$
and $h_1(x):=1$, where $i=\sqrt{-1}$. We recall that $p_i(x)$
corresponds to the coefficient of degree $i$ in $z$ of the
polynomial defining $\dot{x}$ and $q_i(x)$ corresponds to the
coefficient of degree $i$ in $z$ of the polynomial defining
$\dot{z}$. We have that $p_4(x) \equiv 0$, so the first assertion
of Proposition \ref{prop15} is satisfied. We compute
$h_1'(x)/h_1(x)+ h_2'(x)$ which gives $2/x^3$ and the values of
the $\tilde{\sigma}_{\kappa}(x)$ are equal to:
$\tilde{\sigma}_0(x)= 0$, $\tilde{\sigma}_1(x)= -2/x^2$,
$\tilde{\sigma}_2(x)= -2/x^4$, $\tilde{\sigma}_3(x)= 6/x^4$,
$\tilde{\sigma}_4(x)= 4/x^6$. An easy computation shows that the
formulas written in (\ref{eq11a}) give $k_0(x)=k_1(x)=k_3(x)
\equiv 0$ and $k_2(x)=-4x$, which corresponds to the given value
of the cofactor $k(x,z)=-4 x z^2$.

\subsection{Darboux functions obtained from invariants\label{sect34}}

We consider a Darboux function $I(x,y)$ which is the product of
invariant algebraic curves up to complex numbers. Then, by
computing $y$-roots, it can be expressed in the form $I(x,y)=h(x)
\prod_{i=1}^{\ell} (y-g_i(x))^{\alpha_i}$, where $\alpha_i \in
\mathbb{C}$, $g_i(x)$ are algebraic particular solutions and
$h(x)$ is such that its logarithmic derivative is a rational
function. We recall that the {\em logarithmic derivative} of a
function $h(x)$ is the quotient $h'(x)/h(x)$.
\par
In this subsection, we are concerned with the converse of this
problem, that is, we give the conditions that a function of the
form $I(x,y)=h(x) \prod_{i=1}^{\ell} (y-g_i(x))^{\alpha_i}$ must
satisfy in order to be a Darboux function.

\begin{proposition}
Assume that the function $I(x,y)=h(x) \prod_{i=1}^{\ell}
(y-g_i(x))^{\alpha_i}$, where $\alpha_i \in \mathbb{C}$ and
$g_i(x)$ are algebraic functions, and $h(x)$ is such that
$h'(x)/h(x)$ is a rational function, satisfies that
\[ P(x,y) \frac{\partial I}{\partial x}(x,y) + Q(x,y) \frac{\partial
I}{\partial y}(x,y) = k(x,y) \, I(x,y), \] where $k(x,y)$ is a
polynomial, then $I(x,y)$ is a Darboux function. \label{prop16}
\end{proposition}
{\em Proof.} By Proposition $7$ in \cite{GaGG}, we have that each
$g_i(x)$ is a particular solution of equation (\ref{eq1}), which
is algebraic by hypothesis. Therefore, by Theorem \ref{th6}, we
deduce that the irreducible polynomial $f(x,y)$ such that
$f(x,g(x)) \equiv 0$ gives rise to an invariant algebraic curve of
system (\ref{eq1}). For each $g_i(x)$ there is an invariant
algebraic curve $f(x,y)=0$, but each invariant algebraic curve
$f(x,y)=0$ can implicitly define several functions $g_i(x)$, as
much as the degree $s$ of $f(x,y)$ in $y$. We denote by
$g_{j_i}(x)$ all the algebraic functions defined by the same
invariant algebraic curve $f_j(x,y)=0$. Assume that the finite set
$\{g_i(x) \, : \, i=1,\ldots, \ell\}$ defines a total of $r$
invariant algebraic curves $f_j(x,y)=0$, $j=1,2, \ldots, r$.
Hence, the function $I(x,y)$ is written as $I(x,y)= h(x) \prod
(y-g_{j_i}(x))^{\alpha_{j_i}}$ where the product is taken over all
the possible subindexes. Some of the $\alpha_{j_i}$ can be null in
this denomination.
\par
By Proposition \ref{prop8}, we have that each particular solution
$y=g_{j_i}(x)$ has an associated quasipolynomial cofactor, which
will be denoted by $M_{j_i}(x,y)$. Each invariant algebraic curve
$f_j(x,y)=0$ has a polynomial cofactor denoted by $k_{f_j} (x,y)$.
We define $M_0(x,y):=h'(x) P(x,y)/h(x)$, which is a rational
function by hypothesis. Each polynomial $f_j(x,y)$ decomposes, by
Theorem \ref{th5}, as $f_j(x,y)= h_j(x) \prod_{i=1}^{s_j}
(y-g_{j_i}(x))$, where $h_j(x)$ is a function of $x$ such that
$h_j'(x)/h_j(x)$ is a rational function, $g_{j_i}(x)$ are
algebraic functions and $s_j$ is the degree of $f_j(x,y)$ in $y$.
We define $M_{j_0}(x,y):=h_j'(x) P(x,y)/h_j(x)$ which is a
rational function. We have that:
\begin{eqnarray*}
P(x,y) \frac{\partial f_j}{\partial x} + Q(x,y) \frac{\partial
f_j}{\partial y} & = & P(x,y) h_j'(x) \prod_{i=1}^{s_j}
(y-g_{j_i}(x)) + h_j(x) \\ & &  \left[  \sum_{i=1}^{s_j}
\left(-P(x,y) g_{j_i}'(x) + Q(x,y) \right) \prod_{\nu=1, \nu \neq
i}^{s_j} (y-g_{j_\nu}(x)) \right]
\\ & = & \left( M_{j_0}(x,y) + \sum_{i=1}^{s_j} M_{j_i}(x,y)
\right) h_j(x) \prod_{i=1}^{s_j} (y-g_{j_i}(x)),
\\ & = & \left( M_{j_0}(x,y) + \sum_{i=1}^{s_j} M_{j_i}(x,y)
\right) \, f_{j}(x,y).
\end{eqnarray*}
On the other hand,
\[ P(x,y) \frac{\partial f_j}{\partial x}(x,y) + Q(x,y) \frac{\partial
f_j}{\partial y}(x,y) \ =\ k_{f_j} (x,y) \, f_j(x,y), \] for being
$f_j(x,y)=0$ an invariant algebraic curve with cofactor
$k_{f_j}(x,y)$. Therefore, we deduce the following identities:
\begin{equation}
M_{j_0}(x,y) \,  + \,  \sum_{i=1}^{s_j} M_{j_i}(x,y) \, =\,
k_{f_j}(x,y), \label{eq12}
\end{equation}
for $j=1,2,\ldots, r$. We have that each $M_{j_i}(x,y)$ is a
rational function of $x$, a polynomial in $g_{j_i}(x)$ and a
polynomial in $y$, by Proposition $7$ in \cite{GaGG}. We can
change the powers of $g_{j_i}(x)$ which are equal to or higher
than $s_j$ to a linear combination of lower powers by using the
expression $f_j(x,g_{j_i}(x)) \equiv 0$. Therefore, equating the
powers of $y$ (after changing all the powers $g_{j_i}(x)$ to
combinations of $g_{j_i}(x)^\nu$, $0 \leq \nu < s_j$), we get that
the $r$ identities (\ref{eq12}) give a total of  $r (m-1)$ linear
combinations of the powers $g_{j_i}(x)^\nu$, $0 \leq \nu < s_j$,
where $m$ is the degree of system (\ref{eq1}) in the variable $y$.
\par
We consider the identity
\[ P(x,y) \frac{\partial I}{\partial x}(x,y) + Q(x,y) \frac{\partial
I}{\partial y}(x,y) = k(x,y) \, I(x,y), \] $ $from which we deduce
that
\begin{equation}
M_0 + \sum_{i=1}^{s_1} \alpha_{1_i} M_{1_i} + \sum_{i=1}^{s_2}
\alpha_{2_i} M_{2_i} + \ldots + \sum_{i=1}^{s_r} \alpha_{r_i}
M_{r_i} = k. \label{eq13}
\end{equation}
As before, we can change each power $g_{j_i}(x)^\nu$ with $\nu
\geq s_j$ to a linear combination of the powers $g_{j_i}(x)^\nu$,
$0 \leq \nu < s_j$. And from (\ref{eq13}) we deduce an identity
which is a linear combination of the powers $g_{j_i}(x)^\nu$,
$i=1,2, \ldots, s_j$, $0 \leq \nu < s_j$ and $j=1,2, \ldots, r$.
\par
We notice that if $\alpha_{j_i}=\beta_j$, where $\beta_j \in
\mathbb{C}$, for all $1 \leq i \leq s_j$, and $1 \leq j \leq r$,
then the relations given by (\ref{eq12}) make the identity
(\ref{eq13}) compatible with the fact that $k(x,y)$ is a
polynomial in both variables $x$ and $y$. Assume that this is not
the case. Assume that we have $\alpha_{j_\nu} \neq
\alpha_{j_\upsilon}$ for certain $\nu, \upsilon$. We can assume
that $j=1$ without loss of generality. We can consider the $s_j$
symmetric polynomials defined by each $f_j(x,y)$ which are linear
combinations of the powers $g_{j_i}(x)^\nu$, $i=1,2, \ldots, s_j$,
$0 \leq \nu < s_j$, $j=1,2, \ldots, r$. By using the elimination
theory, we can eliminate all the appearances of $g_{j_i}(x)$ with
$j > 1$ in (\ref{eq13}). We obtain in this way a relation only
involving $g_{1_i}(x)$. By using the symmetric polynomials
associated to $f_1(x,y)$, we eliminate all the $g_{1_i}(x)$,
except one, which may be $g_{1_1}(x)$. The resulting relation
$R(x,g_{1_1}(x)) \equiv 0$ can only be of two forms: either
$R(x,y)$ is a multiple of $f_1(x,y)$ or it is not. In the first
case, we have that the relation given by (\ref{eq13}) is a
combination of the symmetric polynomials, which are symmetric with
respect to $g_{1_i}(x)$, $i=1,2,\ldots,s_1$. This symmetry implies
that $\alpha_{j_\nu} = \alpha_{j_\upsilon}$, for all $\nu$ and
$\upsilon$. In the second case, we would get that $R(x,y)$ is a
polynomial such that $R(x,g_{1_1}(x)) \equiv 0$ and the degree of
$R(x,y)$ in $y$ is lower than $s_1$. We recall that we have
already substituted all the appearances of powers of $g_{1_1}(x)$
of higher degree by the corresponding expression given by the
equation $f_1(x,g_{1_1}(x)) \equiv 0$. The existence of a
polynomial $R(x,y)$ such that $R(x,g_{1_1}(x)) \equiv 0$ and the
degree of $R(x,y)$ in $y$ being lower than $s_1$ is a
contradiction with the fact that $f_1(x,y)$ is the only
irreducible polynomial satisfying $f_1(x,g_{1_1}(x)) \equiv 0$,
modulus associates.
\par
Hence, we conclude that the only possibility is that
$\alpha_{j_i}=\beta_j$, where $\beta_j \in \mathbb{C}$, for all $1
\leq i \leq s_j$, and $1 \leq j \leq r$. We have that: \[ I(x,y)=
h(x) \frac{ \prod_{j=1}^{r} f_j(x,y)^{\beta_{j}}}{ \prod_{j=1}^{r}
h_j(x)^{\beta_{j}}} \] and we define $\tilde{h}(x)=h(x)
\prod_{j=1}^{r} h_j(x)^{- \beta_j}$. We notice that since $h(x)$
and $h_j(x)$, $j=1,2, \ldots, r$ satisfy that its logarithmic
derivative is a rational function, $\tilde{h}(x)$ also satisfies
that $\tilde{h}'(x)/\tilde{h}(x)$ is a rational function. By
integration, we obtain that $\tilde{h}(x)$ is a Darboux function.
We deduce that $I(x,y)$ is equal to $\tilde{h}(x) \prod_{j=1}^{r}
f_j(x,y)^{\beta_{j}}$ which is a Darboux function, as we wanted to
show. \bbox
\newline

{\bf Example.} We are going to describe an example of Proposition
\ref{prop16} so as to make the proof clearer. Let us consider the
following planar polynomial differential system:
\begin{equation}
\begin{array}{lll}
\displaystyle \dot{x} & = & \displaystyle -5-5 x +15 y^2 - 6 x^2 y
+ 14 x y^2 - 9 x y^4, \\ \displaystyle \dot{y} & = & \displaystyle
5 + 2 x - 3 y - 2 x y^2 + 6 y^3 - 3 y^5.
\end{array}
\label{eq14}
\end{equation}
This system exhibits two invariant algebraic curves of degree $3$:
$f_1(x,y)=0$ with $f_1(x,y): = y^3-y-x$ and $f_2(x,y)=0$ with
$f_2(x,y):= x y^2-x-1$. Their cofactors are, respectively,
$k_{f_1}(x,y)=-3 (1+2xy - 4 y^2 + 3 y^4)$ and $k_{f_2}(x,y) = 5 \,
k_1(x,y)/3$. We factorize the polynomial $f_1(x,y)$ as $h_1(x)
(y-g_{1_1}(x)) (y-g_{1_2}(x)) (y-g_{1_3}(x))$ where $h_1(x):=1$
and $g_{1_i}(x)$, $i=1,2,3$, are the corresponding $y$-roots of
$f_1(x,y)$. It is easy to see that the polydromy order of
$g_{1_i}(x)$ is $1$. We perform the same computations for
$f_2(x,y)$ and we have that $f_2(x,y) = h_2(x) (y-g_{2_1}(x))
(y-g_{2_2}(x))$ where $h_2(x):=x$, $g_{2_1}(x):=\sqrt{1+1/x}$ and
$g_{2_2}(x) := -\sqrt{1+1/x}$. It is easy to see that the
polydromy order of $g_{2_i}(x)$, $i=1,2$, is $2$. We have that
$y-g_{j_i}(x)$ are algebraic particular solutions and we can
compute their corresponding quasipolynomial cofactors
$M_{j_i}(x,y)$ which are:
\begin{eqnarray*}
M_{1_i}(x,y) & := & \frac{3}{1-3g_{1_i}^{2}(x)} \left( -1 + 5 y +
2 y^2- y^4 + ( 5 + 2 x- 2y + 2 y^3) g_{1_i}(x) + \right. \\ & &
\left. + (1 + 2 x y - 4 y^2 + 3 y^4) g_{1_i}^2(x) \right), \\
M_{2_i}(x,y) & := & \frac{1}{2 x g_{2_i}(x)} \left( 5 y - 10 x - 4
x^2 + 6 x y - 15 y^3 - 6 x y^3 + \right. \\ & & \left. + (5 -15
y^2 - 4 x^2 y + 6 x y^2 - 6 x y^4) g_{2_i}(x) \right) .
\end{eqnarray*}
We have that $g_{1_i}(x)$, $i=1,2,3$ are the $y$-roots of
$f_1(x,y)$ and, hence, they satisfy the following relationships,
given by the symmetric polynomials on the $y$--roots:
\[ \begin{array}{llll} \vspace{0.2cm} \displaystyle V_{1_1} & := & g_{1_1}(x) + g_{1_2}(x) +
g_{1_3}(x) & = 0, \\ \vspace{0.2cm} \displaystyle V_{1_2} & := &
g_{1_1}(x) g_{1_2}(x) + g_{1_1}(x) g_{1_3}(x) + g_{1_2}(x)
g_{1_3}(x) & = -1, \\ \displaystyle V_{1_3} & := & g_{1_1}(x)
g_{1_2}(x) g_{1_3}(x) & = x.
\end{array} \]
In the same way we have that $V_{2_1}:= g_{2_1}(x) + g_{2_2}(x) =
0$ and $V_{2_2}:= g_{2_1}(x) g_{2_2}(x) = - 1 -1/x$. These
relationships give that $M_{1_1} + M_{1_2} + M_{1_3} = k_1$ and
$M_{2_0} + M_{2_1}+ M_{2_2}= k_2$ where $M_{2_0}=h_2'(x) \dot{x} /
h_2(x) = (-5-5 x +15 y^2 - 6 x^2 y + 14 x y^2 - 9 x y^4) / x$. We
consider a function of the form: \[ I(x,y) = h(x) (y -
g_{1_1}(x))^{\alpha_{1_1}} (y - g_{1_2}(x))^{\alpha_{1_2}} (y -
g_{1_3}(x))^{\alpha_{1_3}} (y - g_{2_1}(x))^{\alpha_{2_1}} (y -
g_{2_2}(x))^{\alpha_{2_2}}, \] where $h(x)$ is such that
$h'(x)/h(x)$ is a rational function (we define $M_0(x,y) := h'(x)
\dot{x} /h(x)$), $\alpha_{j_i} \in \mathbb{C}$ and $g_{j_i}(x)$
are the functions defined above. We also assume that
$\mathcal{X}\left( I(x,y) \right) = k(x,y) I(x,y)$ where $k(x,y)$
is a polynomial in $x$ and $y$ of degree at most $4$. We have
that:
\[ M_0 + \alpha_{1_1} M_{1_1} + \alpha_{1_2} M_{1_2} +
\alpha_{1_3} M_{1_3} + \alpha_{2_1} M_{2_1} + \alpha_{2_2} M_{2_2}
= k. \] This identity gives rise to five equations relating the
$g_{j_i}(x)$ and $\alpha_{j_i}$ when equating the coefficients of
different degrees of $y$. The equation corresponding to the degree
$4$ in $y$ is equal to: \begin{equation} -3\left( \alpha_{1_1}+
\alpha_{1_2}+ \alpha_{1_3}+ \alpha_{2_1} + \alpha_{2_2} + 3 x \,
\frac{h'(x)}{h(x)}\right) = k_4(x), \label{eqexy4}
\end{equation} where $k_j(x)$ is the coefficient of degree $j$
in $y$ of $k(x,y)$. We note that $k_4(x)$ is a real number since
$k(x,y)$ is a polynomial in $x$ and $y$ of degree at most $4$.
Then, equation (\ref{eqexy4}) implies that $h(x)=x^b$ with $b=
-(\alpha_{1_1}+ \alpha_{1_2}+ \alpha_{1_3}+ \alpha_{2_1} +
\alpha_{2_2})/3 - k_4/9$. The equation corresponding to the degree
$3$ in $y$ gives:
\[ \sum_{i=1}^{3} \left( \frac{6
\alpha_{1_i} g_{1_i}(x)}{1- 3 g_{1_i}(x)^2} \right) - \frac{3 (5 +
2 x)}{2 x} \left( \frac{\alpha_{2_1}}{g_{2_1}(x)} +
\frac{\alpha_{2_2}}{g_{2_2}(x)} \right) = k_3(x).
\]
The following equations correspond to the degrees $2$, $1$ and $0$
in $y$, respectively:
\[
(15+14x) \, \frac{h'(x)}{h(x)} \, + 6 \sum_{i=1}^{3} \left(
\frac{\alpha_{1_i}(2 g_{1_i}(x)^2 -1)}{3 g_{1_i}(x)^2 -1} \right)
+ \frac{3 (2x-5)}{2 x} (\alpha_{2_1} + \alpha_{2_2}) = k_2(x),
\]
\[
\begin{array}{l}
\displaystyle 3 \sum_{i=1}^{3} \left( \alpha_{1_i} \frac{5 - 2
g_{1_i}(x) + 2 g_{1_i}(x)^2}{1- 3 g_{1_i}(x)^2} \right) +
\sum_{i=1}^{2} \left( \alpha_{2_i}
\frac{5 + 6 x - 4 x^2 g_{2_i}(x) }{ 2 x g_{2_i}(x) } \right) = \vspace{0.2cm} \\
\displaystyle \qquad \qquad = 6x^2 \, \frac{h'(x)}{h(x)} + k_1(x),
\end{array}
\]
\[ \begin{array}{l}
\displaystyle 3 \sum_{i=1}^{3} \left( \alpha_{1_i} \frac{-1 + (2 x
+ 5) g_{1_i}(x) + g_{1_i}(x)^2}{1- 3g_{1_i}(x)^2} \right) +
\sum_{i=1}^{2} \left( \alpha_{2_i} \frac{5 g_{2_i}(x) - 10 x - 4
x^2}{ 2 x g_{2_i}}\right) = \vspace{0.2cm} \\
\displaystyle \qquad \qquad = 5 (1 + x ) \, \frac{h'(x)}{h(x)} +
k_0(x). \end{array} \] We consider a common denominator in each
one of these equations and by using elimination theory among the
numerators of these equations and the polynomials $V_{j_i}$, we
deduce that the only possibilities are $\alpha_{1_1} =
\alpha_{1_2} = \alpha_{1_3} =: \beta_1$, $\alpha_{2_1} =
\alpha_{2_2} =: \beta_2$, $h(x) = x^{\beta_2}$ and $k(x,y)=-(3
\beta_1 + 5 \beta_2) (1+2xy - 4 y^2 + 3 y^4)$. Hence,
$I(x,y)=f_1(x,y)^{\beta_1} f_2(x,y)^{\beta_2}$, which is a Darboux
function. \newline

We know that taking $y$--roots any Darboux function can be
expressed in the form $I(x,y) = \exp \left\{ h_2(x)
\prod_{k=1}^{r} (y-a_k(x)) \diagup \prod_{j=1}^{s}
(y-\tilde{g}_j(x)) \right\} h_1(x) \prod_{i=1}^{\ell}
(y-g_i(x))^{\alpha_i}$ where $\alpha_i \in \mathbb{C}$, $g_i(x)$,
$\tilde{g}_j(x)$ and $a_k(x)$ are algebraic functions, $h_1(x)$ is
such that its logarithmic derivative is a rational function and
$h_2(x)$ is a rational function. Next Proposition gives the
converse of this assertion, that is, we give the conditions that a
function of the form $I(x,y)$ must satisfy in order to be a
Darboux function.
\begin{proposition}
Assume that the function
\[ I(x,y) = \exp \left\{ h_2(x)
\frac{\prod_{k=1}^{r} (y-a_k(x))}{\prod_{j=1}^{s}
(y-\tilde{g}_j(x))} \right\} h_1(x) \prod_{i=1}^{\ell}
(y-g_i(x))^{\alpha_i}, \] where $\alpha_i \in \mathbb{C}$,
$g_i(x)$, $\tilde{g}_j(x)$ and $a_k(x)$ are algebraic functions,
$h_1(x)$ is such that its logarithmic derivative is a rational
function and $h_2(x)$ is a rational function, satisfies that:
\[ P(x,y) \frac{\partial I}{\partial x}(x,y) + Q(x,y) \frac{\partial
I}{\partial y}(x,y) = k(x,y) \, I(x,y), \] where $k(x,y)$ is a
polynomial, then $I(x,y)$ is a Darboux function. \label{prop17}
\end{proposition}
{\em Proof.} From Theorem 2 in \cite{GG}, we have that each
$g_i(x)$ is a particular solution of (\ref{eq2}), so it has an
associated quasipolynomial cofactor (see Proposition \ref{prop8})
which we denote by $M_i(x,y)$. We define $\Phi(x,y):= \exp
\{h_2(x) \, A_1(x,y) / A_0(x,y) \}$ where $A_1(x,y)=
\prod_{k=1}^{r} (y-a_k(x))$ and $A_0(x,y)= \prod_{j=1}^{s}
(y-\tilde{g}_j(x))$. Since $\mathcal{X}(I) = k I$, we deduce that
\[ \mathcal{X} \left( \Phi(x,y) \right) = \left( k(x,y)-
\frac{h_1'(x)}{h_1(x)} P(x,y) - \sum_{i=1}^{\ell} \alpha_i
M_i(x,y) \right) \Phi(x,y). \] Therefore, we are under the
hypothesis of Proposition \ref{prop12} and Theorem \ref{th13}.
Hence, we realize that the proof of this assertion goes exactly as
the proof of Proposition \ref{prop16}. \bbox

As a consequence of Theorem \ref{th6} and Propositions
\ref{prop14}, \ref{prop15}, \ref{prop16} and \ref{prop17} we can
establish the following result.

\begin{theorem}
Assume that system {\rm (\ref{eq1})} has a first integral or an
integrating factor of the form
\[ I(x,y) = \exp \left\{ h_2(x)
\frac{\prod_{k=1}^{r} (y-a_k(x))}{\prod_{j=1}^{s}
(y-\tilde{g}_j(x))} \right\} h_1(x) \prod_{i=1}^{\ell}
(y-g_i(x))^{\alpha_i}, \] where $\alpha_i \in \mathbb{C}$,
$g_i(x)$, $\tilde{g}_j(x)$ and $a_k(x)$ are algebraic functions,
$h(x)$ and $h_1(x)$ have a rational logarithmic derivative and
$h_2(x)$ is a rational function. Then, $I(x,y)$ is a Darboux
function. \label{th18}
\end{theorem}
{\em Proof.} We are under the hypothesis of Proposition
\ref{prop16} or \ref{prop17}. In this case, $k(x,y)$ is
identically zero or $k(x,y)$ is minus the divergence of the
system. We deduce that $I(x,y)$ must be a Darboux function. \bbox
\par
We notice that when we apply the method of constructing first
integrals given in \cite{GaGG} and all the $h(x)$ and $g_i(x)$ are
completely determined functions, by Theorem \ref{th18} we have
that this first integral is a Darbouxian function. However, the
reciprocal is not true. As some examples in \cite{GaGG} show, we
may have a system with a Darbouxian first integral, but when we
apply the method described in \cite{GaGG} we get a nonlinear
superposition principle. That is, not all the $g_i(x)$ introduced
in the ansatz $I(x,y)$ are determined. If we choose those
undetermined $g_i(x)$ as algebraic particular solutions, we will
have that the first integral given by the superposition principle
becomes a Darboux function.

\vspace{0.5cm}

{\bf Addresses and e-mails:} \\
$^{\ (1)}$ Lab. de Math\'ematiques
et Physique Th\'eorique. CNRS UMR 6083. \\ Facult\'e des Sciences
et Techniques. Universit\'e de Tours. \\ Parc de Grandmont, 37200
Tours, FRANCE.
\\ {\rm E-mail:} {\tt Hector.Giacomini@lmpt.univ-tours.fr}
\vspace{0.2cm} \\
$^{\ (2)}$ Departament de Matem\`atica. Universitat de Lleida. \\
Avda. Jaume II, 69. 25001 Lleida, SPAIN. \\ {\rm E--mails:} {\tt
gine@eps.udl.es}, {\tt mtgrau@matematica.udl.es}


\begin{thebibliography}{99}

\bibitem{Casas} {\sc E. Casas-Alvero},{\it \ Singularities of
Plane Curves.} London Mathematical Society Lecture Note Series,
{\bf 276}. Cambridge University Press, Cambridge, 2000.

\bibitem {ChGGLl} {\sc J. Chavarriga, H. Giacomini, J. Gin\'e and
J. Llibre},{\it \ Darboux integrability and the inverse
integrating factor.} J. Differential Equations {\bf 194} (2003),
116--139.

\bibitem{Christopher1} {\sc C. Christopher},{\it \ Invariant
algebraic curves and conditions for a centre.} Proc. Roy. Soc.
Edinburgh Sect. A {\bf 124} (1994), 1209--1229.

\bibitem{Christopher2} {\sc C. Christopher}, {\it \ Liouvillian first
integrals of second order polynomial differential systems.}
Electron. J. Differential Equations, Vol. {\bf 1999}(1999), 1--7.

\bibitem{ChLl} {\sc C. Christopher and J. Llibre},{\it \ Integrability
via invariant algebraic curves for planar polynomial differential
systems.} Ann. Differential Equations {\bf 16} (2000), 5--19.

\bibitem{GaGG} {\sc I.A. Garc\'{\i}a, H. Giacomini, and J. Gin\'e},
{\it \ Generalized nonlinear superposition principles for
polynomial planar vector fields.} J. Lie Theory, {\bf 15} (2005),
89--104.

\bibitem{GaG}{\sc I.A. Garc\'{\i}a and J. Gin\'e}, {\it \ Generalized
cofactors and nonlinear superposition principles}, Appl. Math.
Lett. {\bf 16} (2003), 1137--1141.

\bibitem{GG} {\sc H. Giacomini and J. Gin\'e},{\it \ An algorithmic
method to determine integrability for polynomial planar vector
fields}, preprint, Universitat de Lleida, 2004.

\bibitem{GGG} {\sc H. Giacomini, J. Gin\'e and M. Grau}, {\it \
Integrability of planar polynomial differential systems through
linear differential equations}, to appear in Rocky Mountain J.
Math.

\bibitem{Hille} {\sc E. Hille},{\it \ Analytic function theory. Vol. {\rm II}.}
Introductions to Higher Mathematics, Ginn and Co. 1962.

\bibitem{Lang} {\sc S. Lang},{\it \ Algebra.} Revised third edition. Graduate Texts in
Mathematics, {\bf 211}. Springer-Verlag, New York, 2002.

\bibitem{Moulin1} {\sc J. Moulin--Ollagnier}, {\it \ About a conjecture on quadratic vector fields.}
J. Pure Appl. Algebra {\bf 165} (2001), 227--234.

\bibitem{Moulin2} {\sc J. Moulin--Ollagnier}, {\it \ Simple Darboux points of polynomial
planar vector fields.} J. Pure Appl. Algebra {\bf 189} (2004),
247--262.

\bibitem{PrelleSinger} {\sc M.J. Prelle and M.F. Singer}, {\it \ Elementary
first integrals of differential equations.} Trans. Amer. Math.
Soc. {\bf 279} (1983), 215--229.

\bibitem{Schinzel} {\sc A. Schinzel} {\it \ Polynomials with special regard
to reducibility.} Encyclopedia of Mathematics and its
Applications, {\bf 77}. Cambridge University Press, Cambridge,
2000.

\bibitem{Singer} {\sc M.F. Singer}, {\it \ Liouvillian first
integrals of differential equations.} Trans. Amer. Math. Soc. {\bf
333} (1992), 673--688.

\bibitem{Walker} {\sc R.J. Walker},{\it \ Algebraic curves.}
Reprint of the 1950 edition. Springer-Verlag, New York-Heidelberg,
1978.

\end{thebibliography}
\end{document}